\pdfoutput=1
\documentclass[3p,a4paper,12pt]{elsarticle}

\usepackage{etex}

\usepackage{graphics}
\usepackage{graphicx}
\usepackage{epsfig}
\usepackage{pst-grad} 
\usepackage{pst-plot} 
\usepackage{pstricks}
\usepackage{lmodern}
\usepackage{xcolor,rotating,epic,eepic}
\usepackage[permil]{overpic}
\usepackage[english]{babel}
\usepackage{lscape}
\usepackage{array}
\usepackage{multirow}
\usepackage{mathrsfs}
\usepackage{amsmath}
\usepackage{amssymb}
\usepackage{float}
\usepackage{subfig}
\usepackage{placeins}
\usepackage{epstopdf}
\usepackage{bm}
\usepackage{helvet}
\usepackage[eulergreek]{sansmath}
\usepackage{booktabs}
\usepackage[hidelinks]{hyperref}

\renewcommand{\Re}[1]{\operatorname{Re}(#1)}
\renewcommand{\Im}[1]{\operatorname{Im}(#1)}
\newcommand*\centermathcell[1]{\omit\hfil$\displaystyle#1$\hfil\ignorespaces}

\makeatletter
\@tfor\next:=abcdefghijklmnopqrstuvwxyzABCDEFGHIJKLMNOPQRSTUVWXYZ\do{%
  \def\command@factory#1{%
    \expandafter\def\csname V#1\endcsname{\mathbf{#1}}
  }
 \expandafter\command@factory\next
}
\makeatother

\newcommand{\Vphi}{\pmb{\phi}}
\newcommand{\VPsi}{\pmb{\Psi}}
\newcommand{\Vzero}{\pmb{0}}
\newcommand{\Vpsi}{\pmb{\psi}}
\newcommand{\Vxi}{\pmb{\xi}}

\journal{Journal of Sound and Vibration}

\begin{document}
\begin{frontmatter}
\title{A spectral characterization of nonlinear normal modes}
\author[lee]{G.I.~Cirillo\corref{cor1}}
\ead{gcirillo@ulg.ac.be}
\author[lee]{A.~Mauroy}
\ead{a.mauroy@ulg.ac.be}
\author[bem]{L.~Renson}
\ead{l.renson@bristol.ac.uk}
\author[lme]{G.~Kerschen}
\ead{g.kerschen@ulg.ac.be}
\author[ccc]{R.~Sepulchre}
\ead{r.sepulchre@eng.cam.ac.uk}
\cortext[cor1]{Corresponding author}
\address[lee]{Systems and Modeling research group, Department of Electrical Engineering and Computer Science, University of Li\`ege, Li\`ege, Belgium.}
\address[bem]{Department of Engineering Mathematics, University of Bristol, Bristol, UK.}
\address[lme]{Space Structures and Systems Lab (S3L), Department of Aerospace and Mechanical Engineering, University of Li\`ege, Li\`ege, Belgium.}
\address[ccc]{Control Group, Department of Engineering, University of Cambridge, Cambridge, UK.}

\begin{abstract}
This paper explores the relationship that exists between nonlinear normal modes (NNMs) defined as invariant manifolds in phase space and the spectral expansion of the Koopman operator. Specifically, we demonstrate that NNMs correspond to zero level sets of specific eigenfunctions of the Koopman operator. Thanks to this direct connection, a new, global parametrization of the invariant manifolds is established. Unlike the classical parametrization using a pair of state-space variables, this parametrization remains valid whenever the invariant manifold undergoes folding, which extends the computation of NNMs to regimes of greater energy. The proposed ideas are illustrated using a two-degree-of-freedom system with cubic nonlinearity.
\end{abstract}

\end{frontmatter}

\section{Introduction}\label{sec:intro}
Nonlinear normal modes (NNMs) of vibration can be regarded as the extension of the concept of linear normal modes to nonlinear systems. They have proved useful in a number of applications going from model reduction~\cite{touze2006nonlinear} and the study of localization phenomena~\cite{vakakis1996normal} to system identification~\cite{peeters2011modal,Renson16}. NNMs of undamped systems were first defined as a vibration of unison of the system~\cite{rosenberg1960normal}. The NNM definition was later generalized to encompass modal interactions~\cite{lee2005complicated}. Inspired by the center manifold technique, Shaw and Pierre extended the NNM concept to damped systems~\cite{shaw1993normal}. They defined NNMs as two-dimensional invariant manifolds in phase space which are tangent to the planes formed by the linear modes at the equilibrium point. Such manifolds are invariant under the flow (i.e. trajectories initialized in the manifold do not leave the manifold), a property which extends the invariance of linear modes to nonlinear systems.

Operator-theoretic methods have a long and rich tradition in dynamical systems. In particular, a (nonlinear) system can be studied through the spectral properties of the associated (linear) Koopman operator. Koopman operators were first defined for Hamiltonian systems in~\cite{koopman1931hamiltonian}. In the seminal work of Mezi\'{c} \cite{mezic2005spectral} the spectral properties of the operator were connected to geometrical features of ergodic (i.e. conservative) systems. More recently, this spectral approach was extended to the case of dissipative systems in the context of model reduction \cite{mauroy2013isostables} and stability analysis \cite{mauroy2014global}. The Koopman operator framework is also conducive to data analysis (see e.g. Dynamic Mode Decomposition~\cite{Tu14dynamic}) and was successfully used to capture coherent structures (called Koopman modes) or detect instabilities in various fields (e.g. fluid flows \cite{rowley2009spectral}, power grids \cite{susuki2013nonlinear}, energy-efficient buildings \cite{eisenhower2010decomposing}).

Our aim is to explore the connections that exist between NNMs and the Koopman operator framework. In this direction, a fundamental result reported by Mezi\'{c} in~\cite{mezic2005spectral} is that normal modes of linear oscillations have their natural analogs, termed the Koopman modes, in the context of nonlinear dynamics. Specifically, it was shown that, for a hyperbolic fixed point, the dynamics can be linearized in the entire basin of attraction~\cite{lan2013linearization} and that the coordinates for the linearization can be obtained through the eigenfunctions of the Koopman operator~\cite{mezic2013analysis}. This framework therefore seems very appropriate to extend the concept of linear normal modes to nonlinear systems. Indeed, we will demonstrate that there is a direct connection between NNMs and the spectral expansion of the Koopman operator. More precisely, the invariant manifolds are related to the zero level sets of some eigenfunctions of the operator, in direct analogy with the linear case.

Using this new characterization, we will show how issues related to the parametrization of the invariant manifold can be overcome. In the original work of Shaw and Pierre, a single pair of variables (either in physical or modal coordinates) was used to obtain a parametrization of the manifold, which can possibly hold only locally. For instance, the parametrization may fail when the invariant manifold presents a complex geometry with foldings, as shown in~\cite{jiang2005theconstruction,renson2014effective}. In contrast, interpreting the concept of NNMs in the Koopman operator framework allows us to obtain a parametrization that is valid globally.

The paper is organized as follows. Sections~\ref{sec:NNMintro} and~\ref{sec:KO} review the theory of NNMs and of the Koopman operator, respectively. In Section~\ref{sec:NNMs} the connection between the two frameworks is investigated, first in the linear case and then in the nonlinear case. We also derive a partial differential equation describing the geometry of the invariant manifold and an approximation of the manifold using Taylor series. The theoretical developments are illustrated with a two-degree-of-freedom mechanical system with cubic nonlinearity in Section~\ref{sec:exe}. The conclusions of the present study are summarized in Section~\ref{sec:conc}.

\section{Introduction to NNMs}\label{sec:NNMintro}
Consider the first-order equations of motion of an N-degree-of-freedom nonlinear mechanical system
\begin{equation}
\dot{x}_i = y_i, \quad \quad \dot{y}_i = f(\mathbf{x},\mathbf{y}) \quad \quad i = 1, ..., N \label{eq:EOM}
\end{equation}
where $x_i$ represents a generalized coordinate (displacement or rotation), $y_i$ is the corresponding velocity, $\mathbf{x} = \left[x_1,...,x_N\right]$ and $\mathbf{y} = \left[y_1,...,y_N\right]$. From this point on, we assume that system~\eqref{eq:EOM} is damped, $\left(\mathbf{x},\mathbf{y}\right)=(\mathbf{0},\mathbf{0})$ is an equilibrium point, and the system linearized around $(\mathbf{0},\mathbf{0})$ has $n/2$ pairs of complex-conjugate eigenvalues with non-zero real parts (hyperbolic equilibrium point).

In the 1990s, Shaw and Pierre defined a NNM as a two-dimensional manifold in phase space~\cite{shaw1993normal} that is invariant under the flow, i.e. trajectories with initial conditions in the manifold remain in it for all time. As such, the invariance property of linear normal modes was extended to nonlinear systems: a motion that starts on a mode stays on this mode, while the other modes remain quiescent. Inspired by the center manifold technique~\cite{CarrBook}, Shaw and Pierre mathematically described the two-dimensional invariant manifold using a pair of master coordinates, the other state-space coordinates being functionally related to the master pair. Choosing arbitrarily $(x_i,y_i)$ as master coordinates, the remaining coordinates follow the constraint equations
\begin{equation}
x_j = X_j(x_i,y_i), \quad \quad y_j = Y_j(x_i,y_i), \quad \quad j=1,...,N; \; j \neq i.
\label{eq:mastercoord}
\end{equation}
In several other works (see, for instance,~\cite{Pesheck02}), a linear change of coordinates was used to express Eq.~\eqref{eq:mastercoord} in modal space and describe the invariant manifold as a function of a linear modal coordinate and its corresponding modal velocity.

A set of equations for the functions $X_j$ and $Y_j$ can be derived by substituting the time derivative of Eq.~\eqref{eq:mastercoord} into the equations of motion~\eqref{eq:EOM}. This leads to a set of $2N-2$ coupled partial differential equations (PDEs)
\begin{eqnarray}\label{eq:functional}
Y_j(x_i,y_i) &=& \frac{\partial X_j(x_i,y_i)}{\partial x_i} y_i + \frac{\partial X_j(x_i,y_i)}{\partial y_i} f_i, \notag \\
f_j &=& \frac{\partial Y_j(x_i,y_i)}{\partial x_i} y_i + \frac{\partial Y_j(x_i,y_i)}{\partial y_i} f_i,
\end{eqnarray}
with $f_j = f_j(x_i,\mathbf{X},y_i,\mathbf{Y})$, $\mathbf{X} = \left\{X_j(x_i,y_i)\right\}$, $\mathbf{Y} = \left\{Y_j(x_i,y_i)\right\}$, and $j=1,...,N; \; j \neq i$. The PDEs~\eqref{eq:functional} do not depend on time and can be solved for $\mathbf{X}$ and $\mathbf{Y}$. Finally, after a substitution of the solutions in the ordinary differential equations governing the master coordinates, the dynamics on the NNM is reduced to a single-degree-of-freedom (SDOF) nonlinear oscillator dynamics
\begin{equation}\label{eq:nl_modal_dynamics}
\dot{x_i} = y_i, \quad \quad \dot{y_i} = f_i(x_i,\mathbf{X},y_i,\mathbf{Y}).
\end{equation}

\begin{figure}[htp]
\centering
\includegraphics[width=0.6\textwidth]{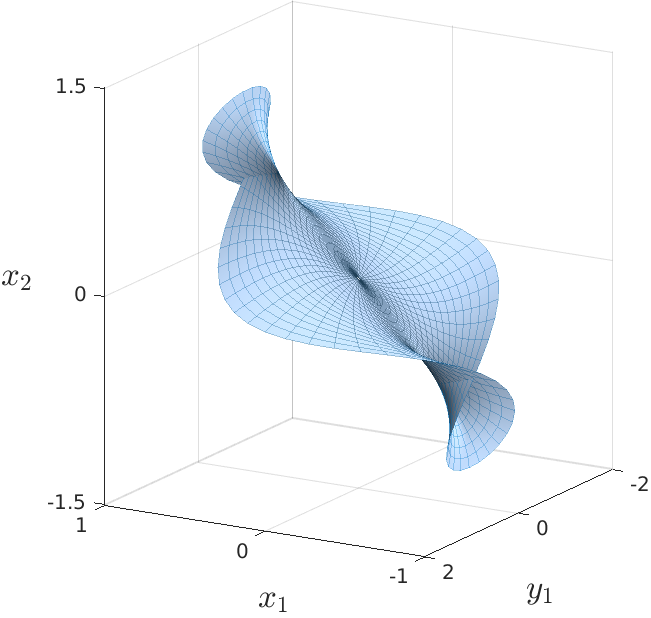}
\caption{In-phase NNM of a two-degree-of-freedom system with cubic nonlinearity. The manifold was computed using the Taylor series expansion presented in Section~\ref{sec:taylor} .}
\label{fig:NNM_def_2DOFMP}
\end{figure}

The PDEs~\eqref{eq:functional} admit several solutions that correspond to the extension of the underlying linear normal modes. The constraint relations $(\mathbf{X},\mathbf{Y})$ describe the geometry of a single NNM in phase space. For a linear system, the mode corresponds to a plane whose dependence in $x_i$ and $y_i$ is function of the inertia, stiffness and damping properties of the system. For a nonlinear system, the mode is no longer flat but curved due to nonlinear distortions. At the equilibrium point, the NNM is tangent to the plane formed by the linear mode of the underlying linear system. Figure~\ref{fig:NNM_def_2DOFMP} illustrates the in-phase NNM of a two-degree-of-freedom system with cubic nonlinearity. The surface shown in the figure corresponds to the position of the second DOF $x_2$ as a function of the master coordinates $(x_1,y_1)$. 
The solutions of~\eqref{eq:functional} can be approximated in the form of a polynomial series expansion in $x_i$ and $y_i$~\cite{shaw1993normal}, whose coefficients are the solutions of a set of algebraic equations. More recently, several numerical methods were proposed in the literature and applied to various conservative~\cite{Pesheck02,Blanc13} and nonconservative~\cite{renson2014effective} systems, see~\cite{RensonReview} for a review. 

Either derived using state space or modal coordinates, Eq.~(\ref{eq:mastercoord}-\ref{eq:nl_modal_dynamics}) are based on the fundamental underlying assumption that the invariant manifold can be described as a function of the chosen pair of variables. However, the invariant manifold may present a complex geometry with several foldings that appear for a number of reasons including internal resonances, localization, and multiple fixed points~\cite{RensonReview,Boivin95ir,Jiang05ir}. Such a manifold cannot generally be described using a simple parametrization that is linear in the state space coordinates. For instance, the invariant manifold shown in Figure~\ref{fig:NNM_def_2DOFMP} exhibits foldings, which prevent a parametrization in terms of the master coordinates $(x_1,y_1)$. In this case, a proper parametrization cannot be obtained as a solution of the PDEs~\eqref{eq:functional}.

In Ref.~\cite{Boivin95ir}, the parametrization issue was addressed using a larger set of master coordinates. When considering $n_m$ pairs, the slave coordinates can be described as
\begin{equation}
x_j = X_j(\mathbf{x}_m,\mathbf{y}_m), \quad \quad y_j = Y_j(\mathbf{x}_m,\mathbf{y}_m),
\label{eq:multi_mode}
\end{equation}
where $(\mathbf{x}_m,\mathbf{y}_m)$ are the vectors of master coordinates. After solving the corresponding manifold-governing PDEs, the dynamics is reduced to $m$ coupled (nonlinear) oscillators. Although this method is elegant and proved effective on a nonlinear beam~\cite{Boivin95ir,Jiang05ir}, it does not completely solve the intrinsic parametrization issue and may therefore fail in some regions of the phase space.

In this paper the parametrization issue is addressed using the relation that exists between NNMs and Koopman operator. Through a nonlinear change of coordinates, the dynamics on the NNM will be expressed as a single-degree-of-freedom linear oscillator. This novel approach bears strong resemblance with the normal form method~\cite{touze2006nonlinear,neild2013nonlinear}, whose main idea is to simplify the vector field using a series of coordinate changes. Unlike the normal form method which has only local validity, the Koopman operator framework provides a change of coordinates that is valid in the entire basin of attraction of a hyperbolic fixed point~\cite{lan2013linearization}.
Therefore, it is possible to think of methods based on Koopman operator as a generalization of methods based on normal form, with the notable difference that the system of coordinates obtained in this case is valid in the whole basin of attraction of a fixed point. By doing so the parametrization issue is overcome, and coupled PDEs describing the NNM in the entire basin of attraction of the equilibrium point can be obtained.

\section{Koopman Operator}\label{sec:KO}
\subsection{Motivation and definition}
\begin{figure}[h]
\centering
\vspace{.5em}
\includegraphics{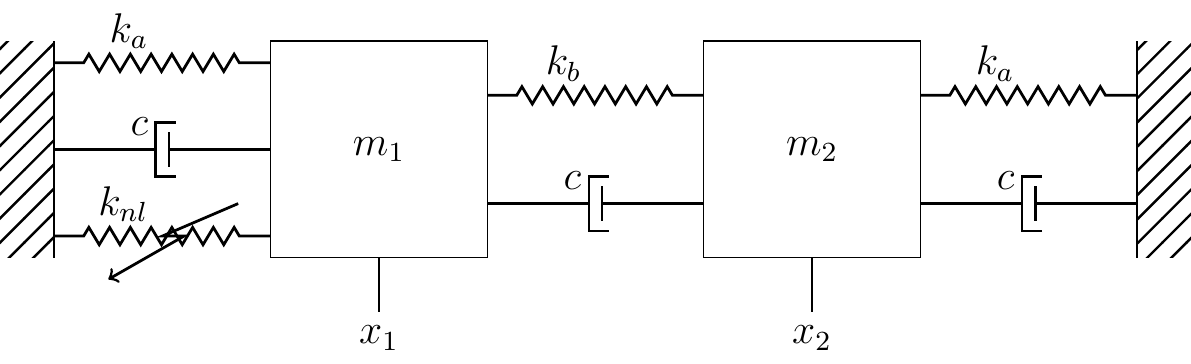}
\vspace{.5em}
\caption{2DOF system with cubic nonlinearity.}\label{fig:SysScheme}
\end{figure}

The system depicted in Figure~\ref{fig:SysScheme} is considered in order to illustrate the Koopman operator definition and properties. The system dynamics can be expressed through a system of ordinary differential equations (ODEs)
\begin{equation}\label{eq:linSysMot}
\dot{\Vx}=\Vf(\Vx),\,\,\,\,\, \Vx=(x_1,x_2,x_3,x_4)^T;
\end{equation}
where $x_1$ and $x_2$ are, respectively, the positions of the first and second mass, and $x_3$ and $x_4$ the corresponding velocities. The general solution, or flow, of~\eqref{eq:linSysMot} is a function, $\Vphi(\Vx_0,t)$, that assigns a point in state space to an initial condition $\Vx_0$ and a time $t$.

It is often interesting to analyze the dynamics of a system in terms of specific quantities such as, for instance, the relative displacement between the first and second mass $g_1(\Vx)=x_2-x_1$, the potential energy $g_2(\Vx)=\frac{1}{2} (k_a(x_1^2+x_2^2)+k_b (x_1-x_2)^2)+\frac{1}{4}k_{nl}x_1^4$ or the kinetic energy $g_3(\Vx)=\frac{1}{2}(m_1x_3^2+m_2x_4^2)$. All these quantities are functions from the state space to $\mathbb{R}$ (or more generally $\mathbb{C}$) and are called \emph{observables}.

Since state-space variables change with time, so do these observables. Given initial conditions $\Vx_0$, the time evolution of an observable is obtained through its composition with the flow. For instance, the potential energy varies according to $V(\Vx_0,t)=g_2(\Vphi(\Vx_0,t))$, while the kinetic energy is given by $K(\Vx_0,t)=g_3(\Vphi(\Vx_0,t))$. At every time instant $t$, $K(\Vx_0,t)$ and $V(\Vx_0,t)$ can still be considered as observables (i.e. functions from the state space to $\mathbb{C}$) such that, generally speaking, the time evolution of an observable results in a family of observables given by its composition with the flow. 

The Koopman operator theory concentrates on how arbitrary observables change with time, instead of how trajectories are organized in state space. Given a dynamical system
\begin{equation}\label{eq:dyn_sys}
\dot{\Vx}=\Vf(\Vx),\,\,\,\,\,\,\Vx\in\mathbb{R}^n,
\end{equation}
and a vector space of observables $\mathcal{F}$, the Koopman operator (or more rigorously, semigroup of operators) $U^t$ is defined as the composition of the flow with $g(\Vx)\in\mathcal{F}$:
\begin{equation}\label{eq:Koop_def}
U^tg(\Vx)=g(\phi(\Vx,t)).
\end{equation}

Equation~\eqref{eq:Koop_def} considers the case of a scalar observable $g$. The definition can be extended to a vector-valued function $\Vg(\Vx)$ by applying $U^t$ to each component $\Vg_i(\Vx)$ of the vector, provided that each vector component $\Vg_i(\Vx)$ is in $\mathcal{F}$.

The fundamental property of Koopman operator is linearity. Recalling the example in Figure~\ref{fig:SysScheme}, the evolution of the total energy of the system, $g_4(\Vx)=g_2(\Vx)+g_3(\Vx)$, is given by the sum of the evolutions of $g_2$ and $g_3$, i.e. $E(\Vx_0,t)=K(\Vx_0,t)+V(\Vx_0,t)$. More generally, applying $U^t$ to the linear combination of two observables $g_\alpha(\Vx)$ and $g_\beta(\Vx)$ that belong to $\mathcal{F}$ yields
\begin{align}
g_\gamma(\Vx)&=a \, g_\alpha(\Vx)+b\, g_\beta(\Vx)\\
\nonumber\\
\begin{split}
U^tg_\gamma(\Vx)&=g_\gamma(\Vphi(\Vx,t))=a \, g_\alpha(\Vphi(\Vx,t)) +b \, g_\beta(\Vphi(\Vx,t))=\\
&=a \, U^tg_\alpha(\Vx)+b \, U^tg_\beta(\Vx).
\end{split}
\end{align}

The linearity property of the Koopman operator is general and does not depend on the dynamical system. This property is very attractive, especially in the context of vibration analysis, where eigenvalues and eigenvectors form the basis of vibration modes. The Koopman operator allows such tools to be used with nonlinear systems, but the price to pay is that infinite dimensional vector spaces have to be considered.

\subsection{Spectral properties}\label{sec:defBasProp}
From this point on the space of observables $\mathcal{F}$ is fixed. Although the only requirement is invariance with respect to the action of $U^t$ (given $g\in\mathcal{F}$, $U^tg\in\mathcal{F}$ for every $t$), different choices of $\mathcal{F}$ result in different spectral properties. In this work, $\mathcal{F}$ is taken as the space of $C^1$ functions that are analytic in a neighborhood of the origin. Moreover, only systems with a stable hyperbolic fixed point (assumed to be the origin) are considered, and our investigations are restricted to the basin of attraction of that point.

Based on the above assumptions, this section presents the Koopman operator spectral properties that are essential to this paper. A more comprehensive introduction can be found in Ref.~\cite{mezic2013analysis}.

\subsubsection{Eigenfunctions}\label{sec:eigenFunct}
Since $U^t$ is linear, it is possible to consider its eigenfunctions: those observables $s(\Vx)\in \mathcal{F}$ whose evolution in time is given by
\begin{equation}  \label{eq:eigFunDef}
U^ts(\Vx)=s(\Vphi(t,\Vx))=\exp(\lambda t)s(\Vx) \,\,\,\, \forall t,
\end{equation}
in which $\lambda \in \mathbb{C}$ is the corresponding eigenvalue.

As an example consider a simple one-dimensional system
\begin{equation}\label{eq:simpleSys}
  \dot{x}=-\eta x,\,\,\,\,\, \eta>0,
\end{equation}
with corresponding flow $\phi(x,t)=x\exp(-\eta t)$. Applying the Koopman operator to the observable $s_1(x)=x$ yields
\begin{equation}
U^ts_1(x)=s_1(\phi(x,t))=x\exp(-\eta t)=\exp(-\eta t)s_1(x),
\end{equation}
thus, the identity function $s_1(x)=x$ is an eigenfunction with eigenvalue $-\eta$.

Now consider a general linear system defined by a matrix $\VA$. The eigenvalues of $\VA$ are also part of the spectrum of the Koopman operator. Moreover, this property extends to nonlinear systems: if $\VA$ is the Jacobian matrix evaluated at the origin of the vector field defined in~\eqref{eq:dyn_sys}, that is
\begin{equation}
\VA=\left.\frac{\partial \Vf}{\partial \Vx}\right|_{\Vx=0},
\end{equation}
then its eigenvalues $\lambda_k$ are part of the spectrum of the Koopman operator~\cite{mauroy2014global}. In fact, if $\VA$ admits $n$ linearly independent eigenvectors, its eigenvalues are associated with $n$ independent eigenfunctions $s_k(\Vx)$ that are defined in the entire basin of attraction of the equilibrium. 

However, the spectrum of $U^t$ is much richer than the spectrum of $\VA$: it always contains an infinity of eigenvalues. For instance, in the case of system~\eqref{eq:simpleSys}, the observable $s_2(x)=x^2$ is also an eigenfunction:
\begin{equation}
U^ts_2(x)=(\phi(x,t))^2=(x\exp(-\eta t))^2=x^2\exp(-2\eta t)=\exp(-2\eta t)s_2(x).
\end{equation}
A similar reasoning can be applied to each integer power $n$ of $x$, so that $-n\eta$ and $s_n(x)=s_1^n(x)=x^n$ are eigenvalues and eigenfunctions of system~\eqref{eq:simpleSys}, respectively.

It turns out that this construction can be applied to general systems. Let $s_1(\Vx)$ and $s_2(\Vx)$ be two eigenfunctions:
\begin{equation}
\begin{split}
U^ts_1(\Vx)&=\exp(\lambda_1 t)s_1(\Vx)\\
U^ts_2(\Vx)&=\exp(\lambda_2 t)s_2(\Vx).
\end{split}
\end{equation}
Applying the Koopman operator to their product $s_{1,2}(\Vx)=s_1(\Vx)s_2(\Vx)$ yields
\begin{equation}\label{eq:multProp}
\begin{split}
U^ts_{1,2}(\Vx)&=s_{1,2}(\Vphi(\Vx,t))=s_1(\Vphi(\Vx,t))s_2(\Vphi(\Vx,t))=\\
&=\exp(\lambda_1 t)s_1(\Vx)\exp(\lambda_2 t)s_2(\Vx)=\exp((\lambda_1+\lambda_2) t)s_1(\Vx)s_2(\Vx).
\end{split}
\end{equation}
Thus, the sum of two eigenvalues is also an eigenvalue, and equivalently the product of two eigenfunctions is an eigenfunction. Using~\eqref{eq:multProp} it is possible to construct an entire system of eigenfunctions.
Considering the basis eigenfunctions $s_k$ associated with the eigenvalues $\lambda_k$ of the Jacobian matrix, we can obtain an infinity of eigenfuctions and eigenvalues:
\begin{equation}\label{eq:high_ord_eig}
s_{k_1,\ldots,k_n}(\Vx)=s_1^{k_1}(\Vx)\ldots s_n^{k_n}(\Vx),\,\,\,\,\,\,\lambda_{k_1,\ldots,k_n}=k_1\lambda_{1}+\ldots+k_n\lambda_{n}.
\end{equation}

In the following a particular property of the eigenfunctions will be central to establish a connection to NNMs defined in Section~\ref{sec:NNMintro}. Let $s(\Vx)$ be an eigenfunction and $\lambda$ the corresponding eigenvalue; then, its zero level set $\{\Vx\in\mathbb{R}^n\mid s(\Vx)=0\}$ is invariant for the dynamics~\eqref{eq:dyn_sys}~\cite{mauroy2014global}. Indeed, given an initial condition $\Vx_0$ satisfying
\begin{equation}  \label{eq:zeroEigFun}
s(\Vx_0)=0,
\end{equation}
the corresponding trajectory satisfies
\begin{equation}  \label{eq:zero_cond_eigf}
s(\Vphi(t,\Vx_0))=\exp(\lambda t)s(\Vx_0)=0.
\end{equation}
If the eigenvalue associated with $s(\Vx)$ is real, then so is $s(\Vx)$, and~\eqref{eq:zeroEigFun} corresponds to one equality. If $\lambda$ has a nonzero imaginary part, $s(\Vx)$ cannot be purely real, so that~\eqref{eq:zeroEigFun} is a complex equation yielding two equalities $\Re{s(\Vx)}=0$ and $\Im{s(\Vx)}=0$. In the latter case, the complex conjugate $\overline{\lambda}$ is also an eigenvalue with eigenfunction $\overline{s}_{k}(\Vx)$, so that~\eqref{eq:zeroEigFun} corresponds to the zero level set of two (complex conjugate) eigenfunctions.

\subsubsection{Spectral expansion}

One of the main tools in the analysis of linear systems is the expansion of the dynamics in terms of linear modes. For a general linear system
\begin{equation}
\dot{\Vx}=\VA\Vx, \,\,\,\,\,\,\Vx\in\mathbb{R}^n,
\end{equation}
with $\VA$ diagonalizable, the flow can be written as
\begin{equation}\label{eq:specExpLin}
\Vphi(\Vx,t)=\sum_{k=1}^ns_k(\Vx)\Vv_k\exp(\lambda_k t),
\end{equation}
where $\Vv_k$ is a (possibly complex) eigenvector of $\VA$, $\lambda_k$ is the corresponding eigenvalue, and $s_k(\Vx)$ is the projection of $\Vx$ onto the eigenspace generated by $\Vv_k$. Note the use of the same notation for the projections of $\Vx$ and the Koopman operator eigenfunctions; these two concepts will be linked in Section~\ref{sec:linCase}.

The Koopman operator theory allows the decomposition~\eqref{eq:specExpLin} to be generalized to nonlinear systems~\cite{mezic2013analysis}. Let $g$ be an observable that belongs to the span of the eigenfunctions~\eqref{eq:high_ord_eig}. The observable can be expressed as a (possibly infinite) sum of eigenfunctions
\begin{equation}\label{eq:spectrDecHyp}
g(\Vx)=\sum_{\substack{ (k_1,\ldots,k_n)\in\mathbb{N}^n \\ k_1+\ldots+k_n>0}} C_{k_1,\ldots,k_n} s_1^{k_1}(\Vx)\ldots s_n^{k_n}(\Vx).
\end{equation}
Applying $U^t$ to~\eqref{eq:spectrDecHyp}, the time evolution of the observable is given by
\begin{equation}  \label{eq:gen_evo}
U^tg(\Vx)=\sum_{\substack{ (k_1,\ldots,k_n)\in\mathbb{N}^n \\ k_1+\ldots+k_n>0}} C_{k_1,\ldots,k_n} s_1^{k_1}(\Vx)\ldots s_n^{k_n}(\Vx) \exp((k_1 \lambda_1+\cdots+k_n \lambda_n) t).
\end{equation}
The coefficients $C_{k_1,\ldots,k_n}$ are the so-called {\em Koopman modes}. Each of them is the projection of the observable $g$ onto the eigenspace generated by the corresponding eigenfunction $s_1^{k_1}\ldots s_n^{k_n}$. If a vector-valued observable $\Vg:\mathbb{R}^n\to\mathbb{C}^d$ is considered, the Koopman modes in~(\ref{eq:spectrDecHyp}--\ref{eq:gen_evo}) become vector-valued modes $\VC_{k_1,\ldots,k_n}$ and the scalar eigenfunctions remain unchanged.

As an example, consider for system~\eqref{eq:simpleSys} the observable $g_1(x)=x+x^2$. In the previous section it was shown that the functions $s_1^k(x)=x^k$ are eigenfunctions. Since $g_1(x)$ is a finite sum of eigenfunctions, we directly find the corresponding Koopman modes: $C_1 = 1$, $C_2=1$, $C_k=0,\,\, k\ge 3$. A less trivial example of spectral decomposition is given by the function $g_2(x)=\exp(x) - 1$. In this case, the Taylor series
\begin{equation}
  \label{eq:KMexp}
    g_2(x)=\exp(x)-1=\sum_{k=1}^\infty \frac{1}{k!}x^k=\sum_{k=1}^\infty \frac{1}{k!}s_1^k(x),
\end{equation}
shows that the Koopman modes are $C_k=\frac{1}{k!}$. In both cases~\eqref{eq:gen_evo} can be used to obtain the time evolution of the observables.

It is clear that~\eqref{eq:gen_evo} bears strong resemblance with the linear mode expansion~\eqref{eq:specExpLin}. To obtain the same decomposition of the flow as in~\eqref{eq:specExpLin}, it is necessary to find an observable whose evolution corresponds to the flow. This is achieved using the identity function $\Vi\Vd (\Vx)=\Vx$.
In a neighborhood of the origin, $\Vi\Vd (\Vx)$ can be decomposed as sum of eigenfunctions~\cite{mauroy2013isostables}
\begin{alignat}{3}
\Vi\Vd(\Vx)&={} & \centermathcell{\Vx}&=\sum_{\substack{ (k_1,\ldots,k_n)\in\mathbb{N}^n \\ k_1+\ldots+k_n>0}}\Vv_{k_1,\ldots,k_n}s_1^{k_1}(\Vx)\ldots s_n^{k_n}(\Vx). & 
\label{eq:id_dec}
\end{alignat}
Using~\eqref{eq:id_dec} together with~\eqref{eq:gen_evo}, a decomposition of the flow in terms of eigenfunctions is obtained
\begin{alignat}{3}
U^t\Vi\Vd(\Vx)=\Vi\Vd(\Vphi(\Vx,t))&={} & \centermathcell{\Vphi(t,\Vx)}&=\sum_{\substack{(k_1,\ldots,k_n)\in\mathbb{N}^n \\ k_1+\ldots+k_n>0}}\Vv_{k_1,\ldots,k_n}s_1^{k_1}(\Vx)\ldots
s_n^{k_n}(\Vx)\exp((k_1\lambda_1+\ldots+k_n\lambda_n)t). & \notag \\
\label{eq:flow_dec}
\end{alignat}
The decomposition~\eqref{eq:flow_dec} generalizes~\eqref{eq:specExpLin} to nonlinear systems. This infinite series can be truncated to a finite number of terms to approximate the system trajectories. The Koopman modes $\Vv_{k_1,\ldots,k_n}$ associated with the identity function can be computed by solving a sequence of linear systems; the details of the calculation are reported in appendix~\ref{ap:KMcomp}. For linear systems, Eq.~\eqref{eq:flow_dec} coincides with the linear normal mode decomposition, as shown in Section~\ref{sec:linCase}.

\FloatBarrier
\section{Koopman Operator and Vibration Modes}\label{sec:NNMs}
\subsection{The linear case}\label{sec:linCase}
Consider the linear system
\begin{equation}
\dot{\Vx}=\VA\Vx, \,\,\, \Vx\in \mathbb{R}^n,
\end{equation}
with $\VA$ diagonalizable, and the corresponding expansion in linear modes given by
\begin{equation}\label{eq:linFlow2}
\Vphi(\Vx,t)=\sum_{k=1}^n s_k(\Vx)\Vv_k\exp(\lambda_k t).
\end{equation}
The terms $s_k(\Vx)$ are the projections onto the corresponding eigenspace; they can be written in the form
\begin{equation}\label{eq:lin_eigFun}
s_k(\Vx)=\Vw_k^*\Vx,
\end{equation}
where ${}^*$ denotes the conjugate transpose and $\Vw_k$ are the eigenvectors of $\VA^*$, normalized so that
\begin{equation}\label{eq:eigVecNorm}
\Vw_k^*\Vv_h=\left\{\begin{array}{lc}1, & k=h\\0, & k\neq h\end{array}\right. .
\end{equation}
Applying the Koopman operator to $s_k(\Vx)$, using the previous property, yields
\begin{equation}
U^ts_k(\Vx)=s_k(\Vphi(\Vx,t))=\Vw_k^*\sum_{i=1}^n s_i(\Vx)\Vv_i\exp(\lambda_i t)=s_k(\Vx)\exp(\lambda_k t).
\end{equation}
Thus, the projections $s_k(\Vx)$ are eigenfunctions of $U^t$ (which explains the notation used for them). It follows that the expansions of the flow through linear modes and Koopman operator eigenfunctions are identical in the linear case. Also, the general decomposition~\eqref{eq:flow_dec} appears as the natural extension of linear normal mode decomposition to nonlinear systems.

The relation between vibration modes and spectral expansion is well-known in the linear case: given a pair of complex conjugate eigenvalues $\lambda_1$ and $\lambda_2$, the motion on the linear mode is obtained when~\eqref{eq:linFlow2} reduces to
\begin{equation}
s_1(\Vx)\Vv_1\exp(\lambda_1 t)+s_2(\Vx)\Vv_2\exp(\lambda_2 t),
\end{equation}
i.e. when the following conditions on the eigenfunctions are verified:
\begin{equation}\label{eq:linearModeEig}
s_k(\Vx)=0,\,\,\,\,\,\,\,\,k=3,\ldots,n.
\end{equation} 
Thus, the plane defining the linear mode is determined by Eq.~\eqref{eq:linearModeEig}. Given a pair of complex conjugate eigenvalues, the corresponding linear mode is the intersection of the zero level set of all the other eigenfunctions. 

In the next two sections~\eqref{eq:linearModeEig} is generalized to nonlinear systems, and it is shown that this generalization is equivalent to define NNMs in terms of invariant manifolds.

\subsection{Nonlinear coordinate transformation}\label{sec:nnChCo}

The properties of Koopman operator in the general case of a nonlinear system with hyperbolic fixed point can be studied through linearization techniques~\cite{lan2013linearization,mauroy2014global}. The results of~\cite{lan2013linearization} are particularly useful, as they allow a local smooth change of coordinates (or diffeomorphism) $\Vh(\Vx)$ to be extended to the entire basin of attraction of the fixed point. Assuming that the eigenvalues are nonresonant (i.e. relations of the type $\sum_i n_i \lambda_i=0$ with $n_i\in\mathbb{N}$ and $\sum_i n_i >0$ are not possible), there exists in a neighborhood of the origin an analytic diffeomorphism of the form~\cite{arnold1988geometrical}
\begin{equation}\label{eq:local_diff}
\Vz=\Vx+\VZ(\Vx), \,\,\,\, \VZ(\Vzero)=\Vzero,\,\,\,\,\left.\frac{\partial \VZ}{\partial \Vx}\right|_{\Vx=\Vzero}=\Vzero,
\end{equation}
such that
\begin{equation}\label{eq:linzcoord}
\dot{\Vz}=\VA\Vz.
\end{equation}
Note that normal form methods aim at approximating the inverse transformation $\Vx=\Vz+\VX(\Vz)$ through its Taylor series.

According to the results and methods given in~\cite{lan2013linearization}, there exists a linearizing global  $C^1$ diffeomorphism $\Vh(\Vx)$, which is analytic when restricted to a neighborhood of the origin and characterizes the eigenfunctions of the Koopman operator. In particular, the eigenfunction corresponding to the eigenvalue $\lambda_k$ of $\VA$ is given by
\begin{equation}\label{eq:eigFunDiff}
s_k(\Vx)=\Vw_k^*\Vh(\Vx),
\end{equation}
where $\Vw_k$ is a left eigenvector of $\VA$ (here and in the following $\VA$ is assumed diagonalizable). Moreover, it follows from~\eqref{eq:local_diff} that
\begin{equation}
\left.\frac{\partial \Vh}{\partial \Vx}\right|_{\Vx=\Vzero}=\VI,
\end{equation}
so that
\begin{equation}\label{eq:eigGrad}
\left.\frac{\partial s_k}{\partial \Vx}\right|_{\Vx=\Vzero}=\Vw_k^*,
\end{equation} 
i.e. the gradient of the eigenfunction at the origin is a left eigenvector of $\VA$. Note that the eigenfunctions $s_k(\Vx)=\Vw_k^*\Vh(\Vx)$ are related to the coordinates $\Vz$. More precisely, they correspond to coordinates $\Vxi$ given by
\begin{equation}
\begin{split}
\Vxi&=\VW \Vz,\\
\VW&=\left(\begin{array}{c}
\Vw_1^* \\ \vdots \\ \Vw_n^*
\end{array}\right),
\end{split}
\end{equation}
i.e. they are linearizing coordinates, in which the dynamics is diagonal.

\subsection{The nonlinear case}\label{sec:nnCase}
This section generalizes the connection between the eigenfunctions of Koopman operator and linear normal modes of vibration to NNMs as defined by Shaw and Pierre \cite{shaw1993normal}. Recall from Section~\ref{sec:NNMintro} that NNM manifolds are characterized by the following properties:
\begin{enumerate}
\item the manifold is invariant;
\item it passes through a stable fixed point (the origin);
\item the plane tangent to it at the origin is a linear mode of the linearized dynamics.
\end{enumerate}

Following the discussion of Section~\ref{sec:linCase}, consider the surface defined by the intersection of $n-2$ zero level sets of eigenfunctions
\begin{equation}  \label{eq:zeroLevSetInt}
s_k(\Vx)=0, \,\,\,\,\, \forall k\neq k_1,k_2, \,\,\,\,\, \lambda_{k_1}=\overline{\lambda_{k_2}}.
\end{equation}
Relying on the properties of the previous sections, we can show that the manifold defined by~\eqref{eq:zeroLevSetInt} verifies the properties 1-3:
\begin{enumerate}
\item As shown in Section~\ref{sec:eigenFunct} the level sets $s_k(\Vx)=0$ are all invariant, therefore the intersection of any number of them is invariant as well.
\item From~\eqref{eq:local_diff} and~\eqref{eq:eigFunDiff} it follows that
\begin{equation}
s_k(\Vzero)=\Vw_k^*\Vh(\Vzero)=0,
\end{equation} 
thus any zero level set contains the origin, and so does their intersection.
\item The tangent space at the origin of~\eqref{eq:zeroLevSetInt} is defined by the linear equations
\begin{equation}
\frac{\partial s_k}{\partial \Vx}(\Vzero)\Vx=0, \,\,\,\,\, \forall k\neq k_1,k_2
\end{equation}
which, using~\eqref{eq:eigGrad}, reduces to
\begin{equation}
\Vw_k^*\Vx=0, \,\,\,\,\, \forall k\neq k_1,k_2,
\end{equation}
and corresponds to the linear mode associated with $\lambda_{k_1}$ and $\lambda_{k_2}$ as shown in Section~\ref{sec:linCase}.
\end{enumerate}
Therefore, the generalization of linear modes to nonlinear systems in the Koopman operator framework corresponds to the definition of NNMs proposed by Shaw and Pierre.

\subsection{NNM approximation using Taylor series}\label{sec:taylor}
Thanks to the connection established in the previous section, it is possible to use methods from the Koopman operator framework to compute NNMs. In this section a Taylor series approximation based on the spectral expansion~\eqref{eq:flow_dec} is proposed. The underlying idea is to use the values of the eigenfunctions as coordinates. This has the advantage that the dynamics in these coordinates is linear as illustrated in Section~\ref{sec:nnChCo}; furthermore, equations~\eqref{eq:zeroLevSetInt} have a particularly simple form in this coordinate system.

Let $\xi_k=s_k(\Vx)$ so that~\eqref{eq:id_dec} becomes 
\begin{equation}  \label{eq:s-1taylor}
\Vi\Vd(\Vx)=\Vx=\sum_{\substack{ (k_1,\ldots,k_n)\in\mathbb{N}^n \\ k_1+\ldots+k_n>0}}%
\Vv_{k_1,\ldots,k_n}\xi_1^{k_1}\ldots \xi_n^{k_n},
\end{equation}
and~\eqref{eq:zeroLevSetInt}
\begin{equation}\label{eq:zeroLevSetIntXi}
\xi_k=s_k(\Vx)=0, \,\,\,\,\, \forall k\neq k_1,k_2, \,\,\,\,\, \lambda_{k_1}=\overline{\lambda_{k_2}}.
\end{equation}

Suppose $k_1=1$ and $k_2=2$. Equation~\eqref{eq:s-1taylor} expresses the transformation between the coordinate system defined by the Koopman operator eigenfunctions $\Vxi$ and the original coordinates $\Vx$, while~\eqref{eq:zeroLevSetIntXi} defines the NNM in the $\Vxi$ coordinates. Combining the two, we can express the points on the NNM as a function of the variables $\xi_1$ and $\xi_2$:
\begin{equation}  \label{eq:NNM_taylor}
\Vx=\sum_{\substack{ (k_1,k_2)\in\mathbb{N}^2 \\ k_1+k_2>0}}%
\Vv_{k_1,k_2}\xi_1^{k_1}\xi_2^{k_2}:=\VPsi(\xi_1,\xi_2),
\end{equation}
where the notation $\Vv_{k_1,k_2}=\Vv_{k_1,k_2,0,\ldots,0}$ is used. Equation~\eqref{eq:NNM_taylor} assigns to any pair of values of the eigenfunctions $s_1(\Vx)$ and $s_2(\Vx)$ a point on the manifold, and each point on the manifold corresponds to a value $\xi_1=\overline{\xi_2}\in\mathbb{C}$. In other words, the manifold is parametrized using all the admissible values for the eigenfunctions. Moreover, it follows from~\eqref{eq:flow_dec} that the time evolution of a trajectory on the manifold is given by
\begin{equation}\label{eq:trajTaylor}
\Vx(t)=\sum_{\substack{ (k_1,k_2)\in\mathbb{N}^2 \\ k_1+k_2>0}}%
\Vv_{k_1,k_2}\xi_1^{k_1}\xi_2^{k_2}\exp((k_1\lambda_1+k_2\lambda_2)t)=\VPsi(\xi_1(t),\xi_2t(t)).
\end{equation}

Equation~\eqref{eq:NNM_taylor} can be truncated to a fixed order to approximate the invariant manifold in a neighborhood of the origin, and correspondingly~\eqref{eq:trajTaylor} to approximate the NNM motion. The vectors $\Vv_{k_1,k_2}$ are the Koopman modes of the identity function and can be computed as discussed in Appendix~\ref{ap:KMcomp}.

To conclude this section some intrinsic limitations of approximating NNMs with the above Taylor series are discussed. First of all, the eigenfunctions might be analytic only in a neighborhood of the origin, so that the amplitude to which NNMs can be computed using~\eqref{eq:NNM_taylor} might also be limited (see also~\cite{mauroy2014global} for a similar issue in approximating eigenfunctions with Taylor series). Even in cases in which the series converges in the whole basin of attraction, the computation of the Koopman modes is more demanding as the order increases, and becomes unfeasible for very high orders.

Due to this the method proposed in this section is not meant to be a ready-to-use algorithm for the computation of NNMs of arbitrary systems, but rather a simple method to compute NNMs for limited amplitudes in small-dimensional systems. Its main purpose is to show that the theory of Koopman operator can lead to computational methods that use the coordinates given by the eigenfunctions. More generic algorithms could be derived using the results of the following section, in which a system of PDEs resembling~\eqref{eq:functional} is derived.

\subsection{Analogy with Shaw and Pierre's manifold governing PDE}\label{sec:PDEsys}

In Section~\ref{sec:NNMintro} a system of PDEs~\eqref{eq:functional} is derived starting from Shaw and Pierre definition of NNMs. A similar system can be derived for $\VPsi(\xi_1,\xi_2)$ defined in~\eqref{eq:NNM_taylor}. To this end, note that the dynamics of $\xi_k$ is by definition
\begin{equation}
\dot{\xi}_k=\lambda_k\xi_k.
\end{equation}
It follows that the time derivative of $\VPsi(\xi_1(t),\xi_2(t))$ is given by 
\begin{align}
\dot{\Vx}(t)&=\frac{\partial \VPsi(\xi_1(t),\xi_2(t))}{\partial \xi_1}\lambda_1\xi_1(t)+\frac{\partial \VPsi(\xi_1(t),\xi_2(t))}{\partial \xi_2}\lambda_2\xi_2(t).\label{eq:trajNNM4PDE}
\end{align}
Substituting~\eqref{eq:trajNNM4PDE} into the equations of motion~\eqref{eq:dyn_sys}, and evaluating at $t=0$ yields
\begin{equation}\label{eq:NNMpdeI}
\Vf(\VPsi)=\frac{\partial \VPsi}{\partial \xi_1}\lambda_1\xi_1+\frac{\partial \VPsi}{\partial \xi_2}\lambda_2\xi_2.
\end{equation}

Equation~\eqref{eq:NNMpdeI} is a set of $n$ PDEs very similar in principle to the one derived in Section~\ref{sec:NNMintro}. The approaches followed here and in Section~\ref{sec:NNMintro} can be both interpreted as a way to solve
\begin{equation}\label{eq:manifoldInPDE}
s_k(\Vx)=0, \,\,\,\, i=3,\ldots,N.
\end{equation}
When using state-space variables as in~\eqref{eq:functional}, Eq.~\eqref{eq:manifoldInPDE} is inverted, but the inversion might be valid only locally. In contrast, a different set of coordinates is used in~\eqref{eq:NNMpdeI} to parametrize the solution, which avoids issues arising when NNMs have a complex geometry (e.g. folding). The two systems of PDEs~\eqref{eq:functional} and~\eqref{eq:manifoldInPDE} have a similar structure, but~\eqref{eq:functional} corresponds to $n-2$ equations while \eqref{eq:manifoldInPDE} consists of $n$ equalities. Moreover, the second system admits more solutions, due to the fact that the eigenfunctions are defined modulo a scaling factor. Indeed, if $\VPsi(\xi_1,\xi_2)$ is a solution then $\hat{\VPsi}(\xi_1,\xi_2)=\VPsi(c\xi_1,\overline{c}\xi_2)$ is a solution for any complex $c\neq 0$.

To obtain $\VPsi(\xi_1,\xi_2)$ real-valued $\xi_1=\overline{\xi}_2$ must be assumed in~\eqref{eq:NNMpdeI}. Equivalently, a system of equations can be derived using as variables $u=\Re{\xi_1}$ and $v=\Im{\xi_1}$, whose dynamics is
\begin{align}
\dot{u}=\sigma_1 u - \omega_1 v, \\
\dot{v}=\omega_1 u + \sigma_1 v,
\end{align}
with $\lambda_1=\sigma_1+i\omega_1$. This leads to
\begin{equation}\label{eq:NNMpdeR}
\Vf(\Vpsi)=\nabla \Vpsi \left(%
\begin{array}{cc}
\sigma_1 & -\omega_1 \\ 
\omega_1 & \sigma_1%
\end{array}%
\right)\left(%
\begin{array}{c}
u \\ 
v%
\end{array}%
\right),
\end{equation}
with $\Vpsi(u,v)=\VPsi(u+iv,u-iv)$. Clearly, the above procedure can be followed for any pair of complex conjugate eigenfunctions, yielding a system of PDEs~\eqref{eq:NNMpdeI} or~\eqref{eq:NNMpdeR} for each NNM. These systems differ only for the value of $\lambda$ (or $\sigma$ and $\omega$).

\subsection{Uniqueness of NNMs}
\begin{figure}[!b]
\centering
\includegraphics{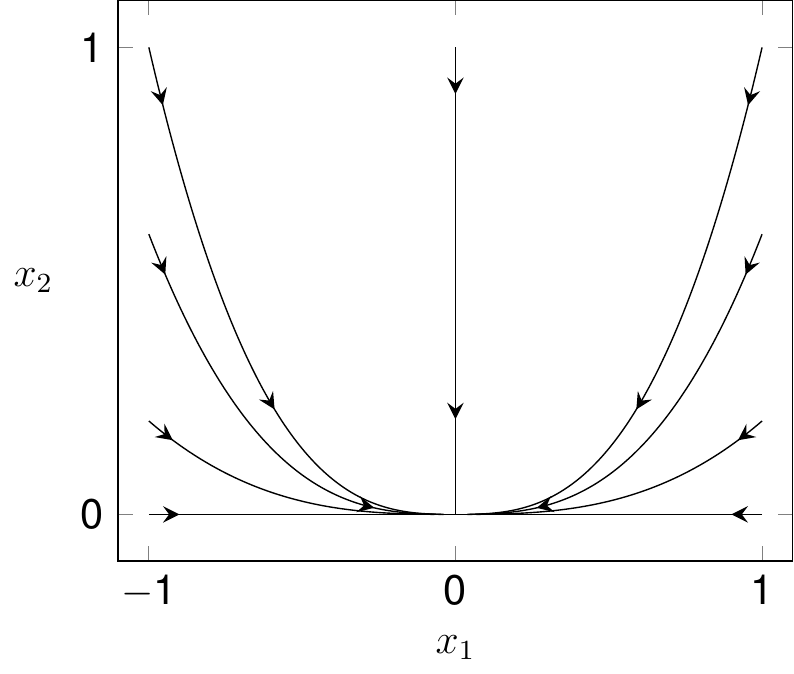}
\caption{Phase portrait of a typical two-dimensional linear system}\label{fig:linPP}
\end{figure}

The section is concluded by a remark about conditions 1-3 of Section~\ref{sec:nnCase} defining the NNMs. Generally these conditions do not determine a unique manifold in the dissipative case. For the sake of illustration, consider a two-dimensional linear system
\begin{equation}
\begin{split}
\dot{x}_1&=\sigma_1 x_1,\\
\dot{x}_2&=\sigma_2 x_2,
\end{split}
\end{equation}
with $\sigma_2<\sigma_1<0$. The phase portrait, illustrated in Figure~\ref{fig:linPP}, shows that almost all trajectories approach the linear eigenspace $x_2=0$, leading to a continuum of tangent manifolds. Analogous phenomena are possible in higher dimensions. For completeness, Appendix \ref{ap:mm} presents a family of invariant manifolds tangent to a linear normal mode, for a four-dimensional linear system. However, assuming nonresonant eigenvalues, the coefficients of a polynomial approximation are uniquely determined at any order, so that only one of the manifolds is analytic. This observation applies to Shaw and Pierre's approach as well (note that the existence of locally analytic manifolds is guaranteed by the existence of a local analytic change of coordinates, see Section~\ref{sec:nnChCo}). The issue of non-uniqueness is discussed in detail in~\cite{Haller16}.

We note that also in the context of the Koopman operator there is the possibility to obtain multiple manifolds for only one eigenvalue. Specifically, our construction of the eigenfunctions relies on the analyticity of the transformation $\Vh$ of Section~\ref{sec:nnChCo}, which is uniquely defined by its Taylor series. However, if this smoothness condition is dropped it is possible to find different transformations verifying~\eqref{eq:local_diff} and~\eqref{eq:linzcoord}, obtaining in this way multiple eigenfunctions for the same eigenvalue.

\section{A 2 DOF example}\label{sec:exe}
The Taylor series expansion derived in Section~\ref{sec:taylor} and based on the Koopman operator framework is used to compute the NNMs of the system shown in Figure~\ref{fig:SysScheme}. The equations of motion are
\begin{equation}\label{eq:sysEq}
\begin{split}
m_1\ddot{x}_1+k_a x_1 +c \dot{x}_1 +k_{nl} x_1^3+k_b(x_1-x_2)+c(\dot{x}_1-\dot{x}_2)&=0,\\
m_2\ddot{x}_2+k_a x_2 +c \dot{x}_2+k_b(x_2-x_1)+c(\dot{x}_2-\dot{x}_1)&=0,
\end{split}
\end{equation}
where $m_1=m_2=m=1$ Kg, $c=0.05$ Ns/m, $k_a=1$ N/m and $k_{nl}=0.5$ N/m$^3$. With $y_1=\dot{x}_1$, $y_2=\dot{x}_2$, $\eta_1=(k_a+k_b)/m$, $\eta_2=k_b/m$, $\beta=c/m$, $\alpha=k_{nl}/m$, \eqref{eq:sysEq} is recast in first-order form
\begin{equation}\label{eq:2doffir}
\left(\begin{array}{c}
\dot{x}_1 \\ \dot{x}_2 \\ \dot{y}_1 \\ \dot{y}_2
\end{array}\right)=\left(\begin{array}{cccc}
0 & 0 & 1 & 0 \\
0 & 0 & 0 & 1 \\
-\eta_1 & \eta_2 & -2\beta& \beta \\
\eta_2 & -\eta_1 & \beta & -2\beta
\end{array}\right)\left(\begin{array}{c}
x_1 \\ x_2 \\ y_1 \\ y_2
\end{array}\right)+\left(\begin{array}{c}
0 \\ 0 \\ -\alpha x_1^3 \\ 0
\end{array}\right).
\end{equation}

\noindent
The Jacobian matrix
\begin{equation}
\VA=\left(\begin{array}{cccc}
0 & 0 & 1 & 0 \\
0 & 0 & 0 & 1 \\
-\eta_1 & \eta_2 & -2\beta& \beta \\
\eta_2 & -\eta_1 & \beta & -2\beta
\end{array}\right)
\end{equation}
has two pairs of complex conjugate eigenvalues $\lambda_1=\overline{\lambda}_2$, $\lambda_3=\overline{\lambda}_4$, which verify a $3:1$ resonance condition $\lambda_1=3\lambda_3$ for $k_b=4$. Exact linear resonance is eliminated by varying $k_b$ but folding is likely to occur in the neighborhood of those parameter values. Table~\ref{tab:eig} shows how the eigenvalues vary for different values of $k_b$.

\begin{table}[h]
\centering
\begin{tabular}{l c c c c c c}
\toprule
$k_b$ & $\lambda_1$ & $f_{n,1}$ & $\zeta_1$ & $\lambda_2$ & $f_{n,2}$ & $\zeta_2$ \\
\midrule
4.7 & -0.075+i3.224 & 3.224 & 0.0233 & -0.025+i0.999 & 1 & 0.025 \\
4.3 & -0.075+i3.079 & 3.098 & 0.0242 & -0.025+i0.999 & 1 & 0.025 \\
4.1 & -0.075+i3.032 & 3.033 & 0.0247 & -0.025+i0.999 & 1 & 0.025 \\
\bottomrule
\end{tabular}
\caption{Eigenvalues, natural frequencies and damping ratios for different values of $k_b$.}\label{tab:eig}
\end{table}

To determine the coefficients describing the NNM in Eq.~\eqref{eq:NNM_taylor}, the trajectory on the manifold given by~\eqref{eq:trajTaylor} is injected into the equations of motion. Evaluating at time $t=0$ yields
\allowdisplaybreaks
\begin{gather}
\sum_{\substack{ (k_1,k_2)\in\mathbb{N}^2 \\ k_1+k_2>0}}\Vv_{k_1,k_2}\xi_1^{k_1}\xi_2^{k_2}(k_1\lambda_1+k_2\lambda_2)=\nonumber\\
\VA\sum_{\substack{ (k_1,k_2)\in\mathbb{N}^2 \\ k_1+k_2>0}}\Vv_{k_1,k_2}\xi_1^{k_1}\xi_2^{k_2}+\left(\begin{array}{c}0 \\ 0 \\ -\alpha\left(\left(\sum\limits_{\substack{ (k_1,k_2)\in\mathbb{N}^2 \\ k_1+k_2>0}}\Vv_{k_1,k_2}\xi_1^{k_1}\xi_2^{k_2}\right)_1\right)^3 \\ 0 \end{array}\right), 
\end{gather} 
where $\left(\Vv\right)_i$ stands for the ith component of the vector $\Vv$. Matching the terms corresponding to the first-order powers results in two systems
\begin{equation}
(\lambda_1\VI-\VA)\Vv_{10}=0,\,\,\,\,\,\,\,\, (\lambda_2\VI-\VA)\Vv_{01}=0,
\end{equation}
so that the first-order modes $\Vv_{10}$ and $\Vv_{01}$ are the eigenvectors of $\VA$ corresponding to $\lambda_1$ and $\lambda_2$. If there is no resonance between the eigenvalues of $\VA$, the successive systems admit a unique solution, and are obtained by matching like-power terms (see the appendix for the complete derivation). Doing so for the second-order terms yields linear equations of the form
\begin{equation}
((k_1\lambda_1+k_2\lambda_2)\VI-\VA)\Vv_{k_1k_2}=0, \,\,\,\, k_1+k_2=2,
\end{equation}
where the right hand side is zero due to the absence of second-order terms in the vector field. We can proceed iteratively, solving for any mode $\Vv_{k_1,k_2}$ a linear system of the form
\begin{equation}
(r\VI-\VA)\Vv_{k_1k_2}=\Vb,\,\,\,\,\,\,\,r=k_1\lambda_1+k_2\lambda_2
\end{equation}
where $\Vb$ is obtained as a sum of lower-order modes (see appendix~\ref{ap:KMcomp} for the details). A change of coordinates using the Jordan basis is used to improve numerical stability.

The above procedure is used to compute the high-frequency (i.e. the out-of-phase) mode of the 2 DOF system for $k_b=4.3$. The Taylor series used to generate the figures of this section were truncated at the $50^{\text{th}}$ order. The maximum amplitude at which the manifolds are represented was chosen so that trajectories obtained with~\eqref{eq:trajTaylor} agree with those computed through numerical integration. To have a quantitative measure of the difference between two trajectories the normalized mean square error (NMSE) was used:
\begin{equation}
\text{NMSE}=\frac{100}{\sigma_x^2}\sum_{k=0}^{K}\|\hat{\Vx}(k)-\Vx(k)\|^2,
\end{equation}
where $\hat{\Vx}$ is the time series obtained with~\eqref{eq:trajTaylor}, $\Vx$ is the reference series obtained through numerical integration, and $\sigma_x^2$ is the variance of the latter. The NMSEs obtained to validate the manifolds were always lower than $1\%$, a value usually taken as threshold for a good match between two time series.

\begin{figure}[!ht]
\centering
\subfloat[]{
\includegraphics[width=.4\textwidth]{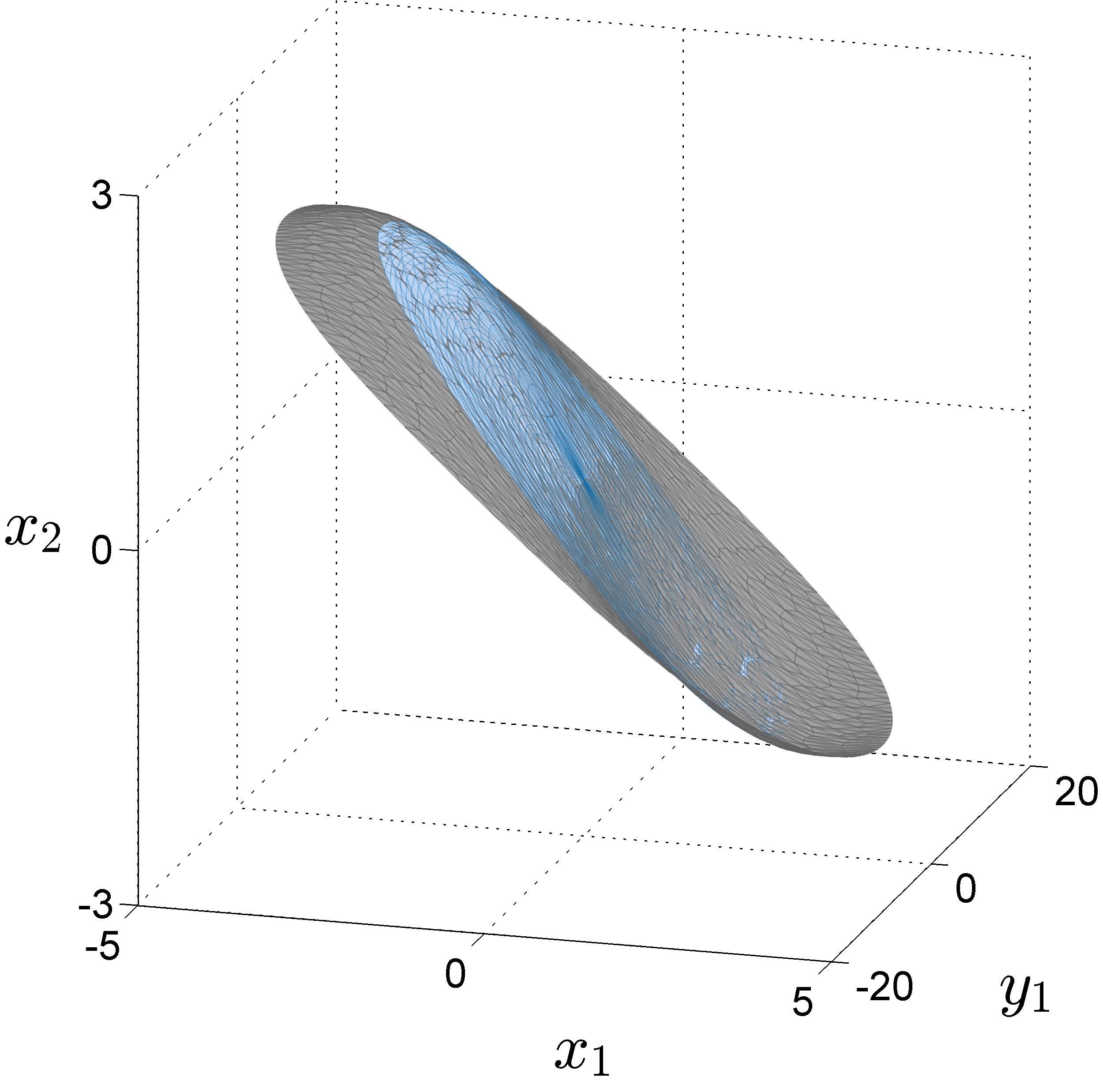}%
}\hfill\subfloat[]{
\includegraphics[width=.37\textwidth]{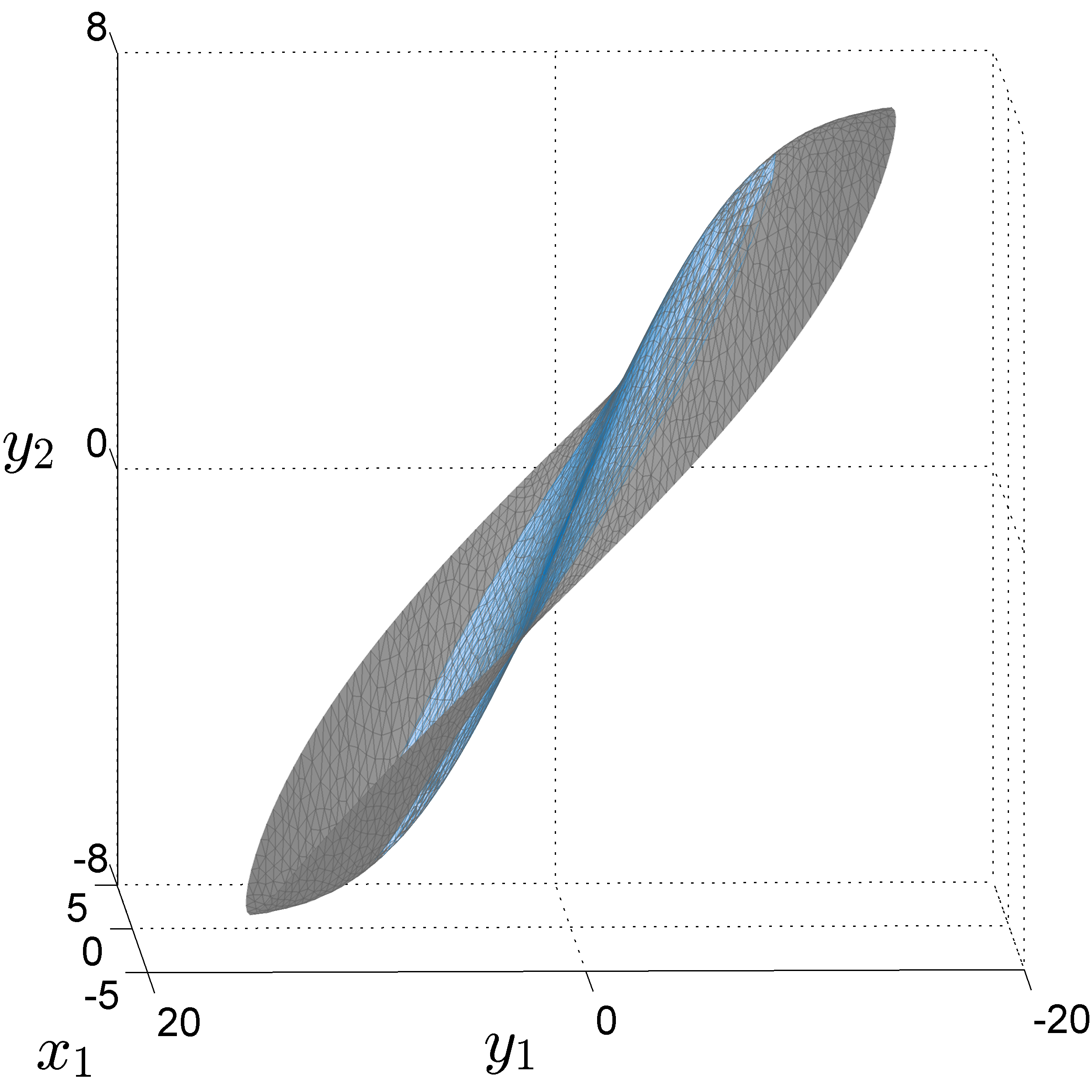}
}
\caption[Out-of-phase NNM]{Out-of-phase NNM;%
\includegraphics{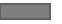}%
FE method,%
\includegraphics{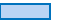}%
Taylor series; (a) $x_1,y_1,x_2$; (b) $x_1,y_1,y_2$}\label{fig:firstNNM}
\end{figure}

\begin{figure*}[p]
\centering
\subfloat[\label{fig:invMan47-1}]{
\includegraphics[height=.28\textheight]{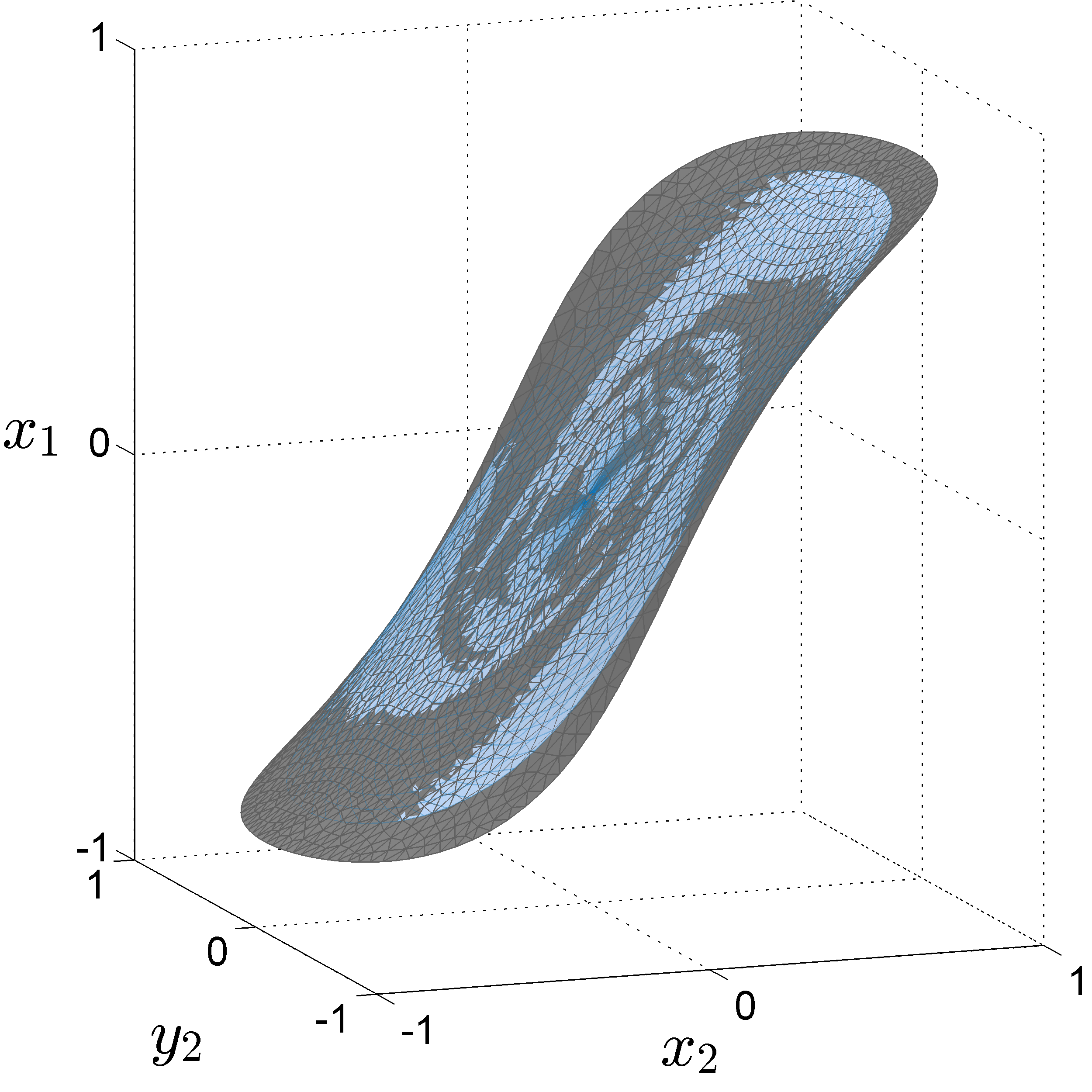}
}\hfill\subfloat[\label{fig:invMan47-2}]{
\includegraphics[height=.28\textheight]{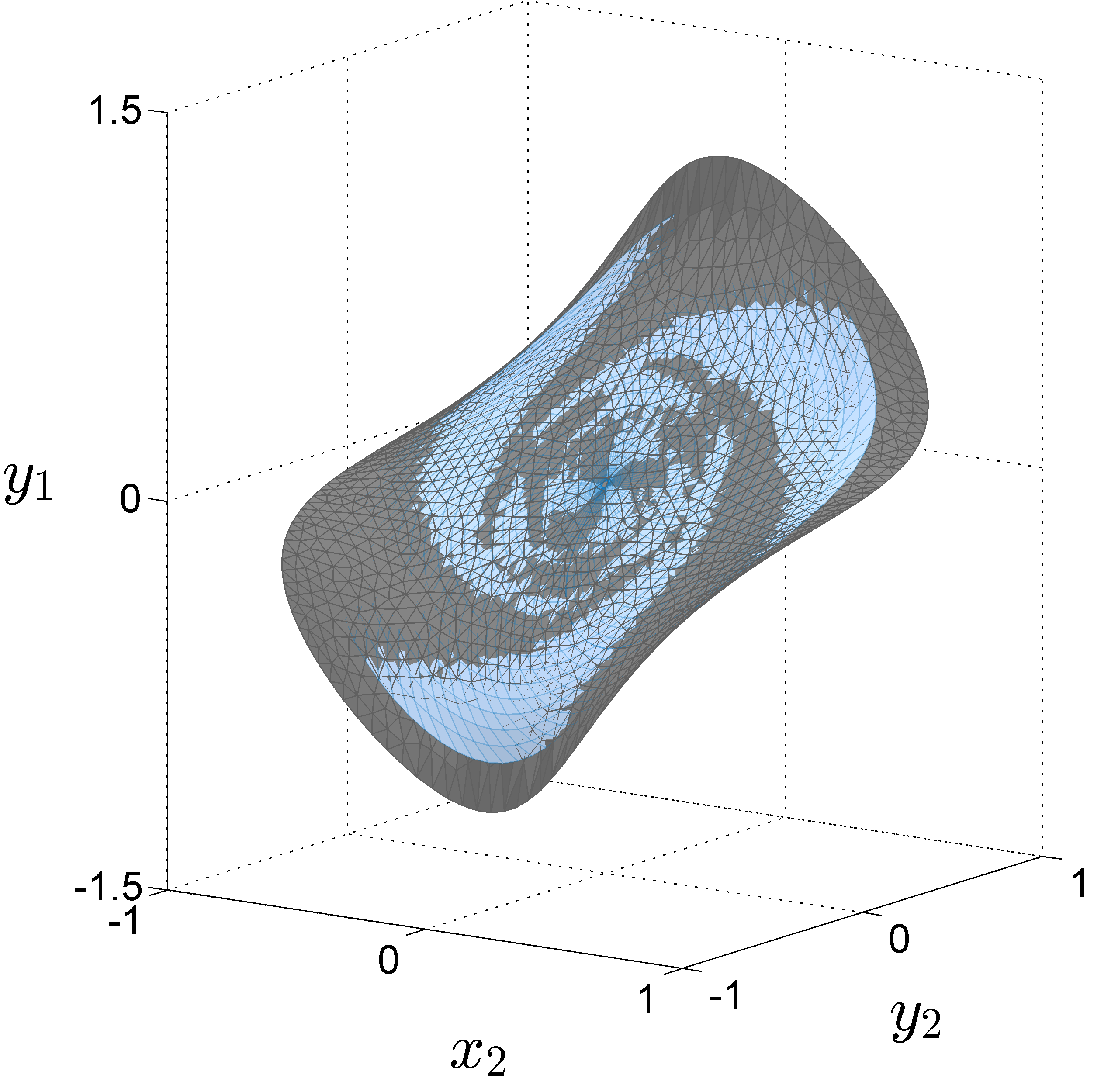}
}

\vfill

\subfloat[\label{fig:invMan43-1}]{
\includegraphics[height=.28\textheight]{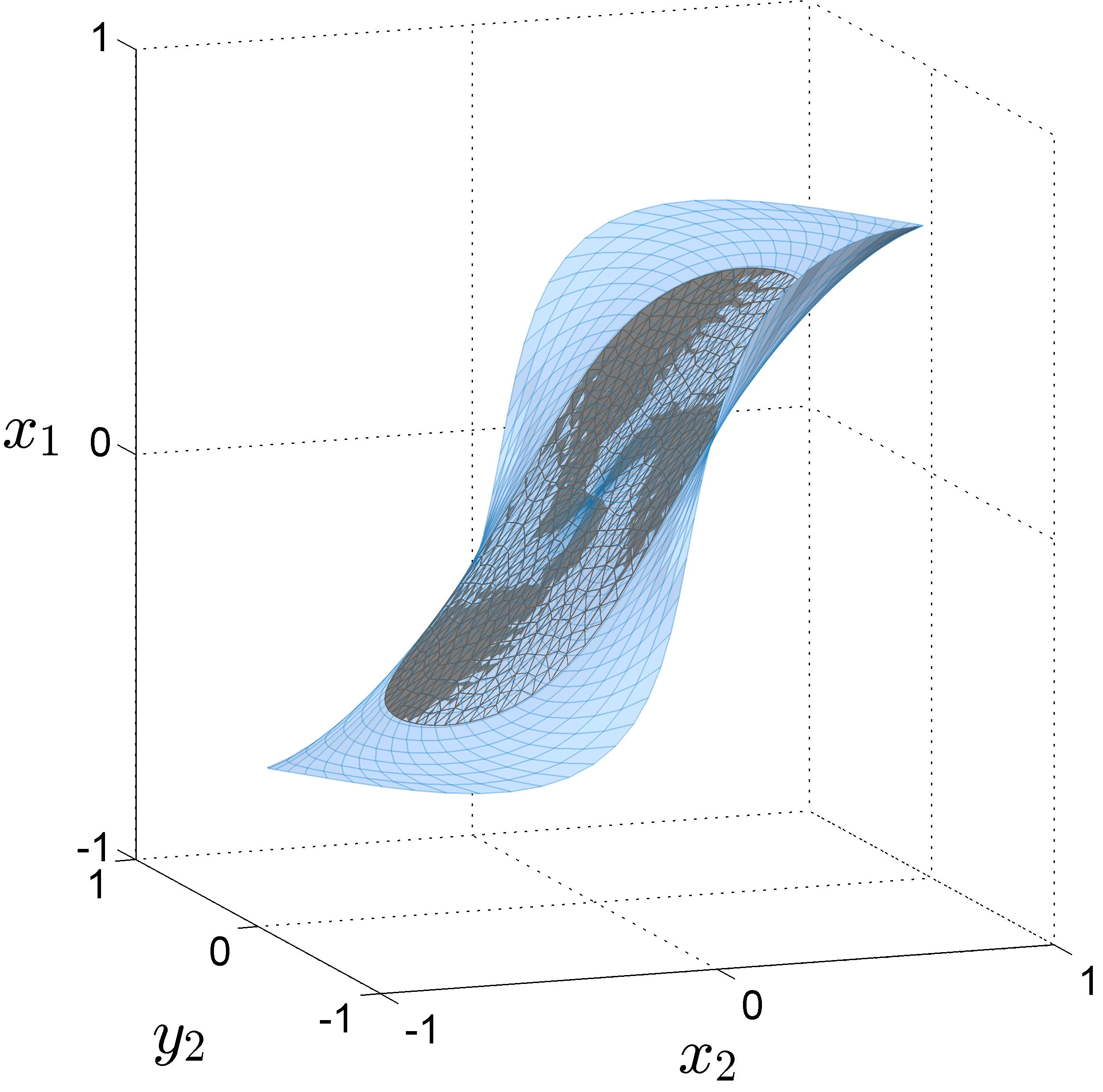}
}\hfill\subfloat[\label{fig:invMan43-2}]{
\includegraphics[height=.28\textheight]{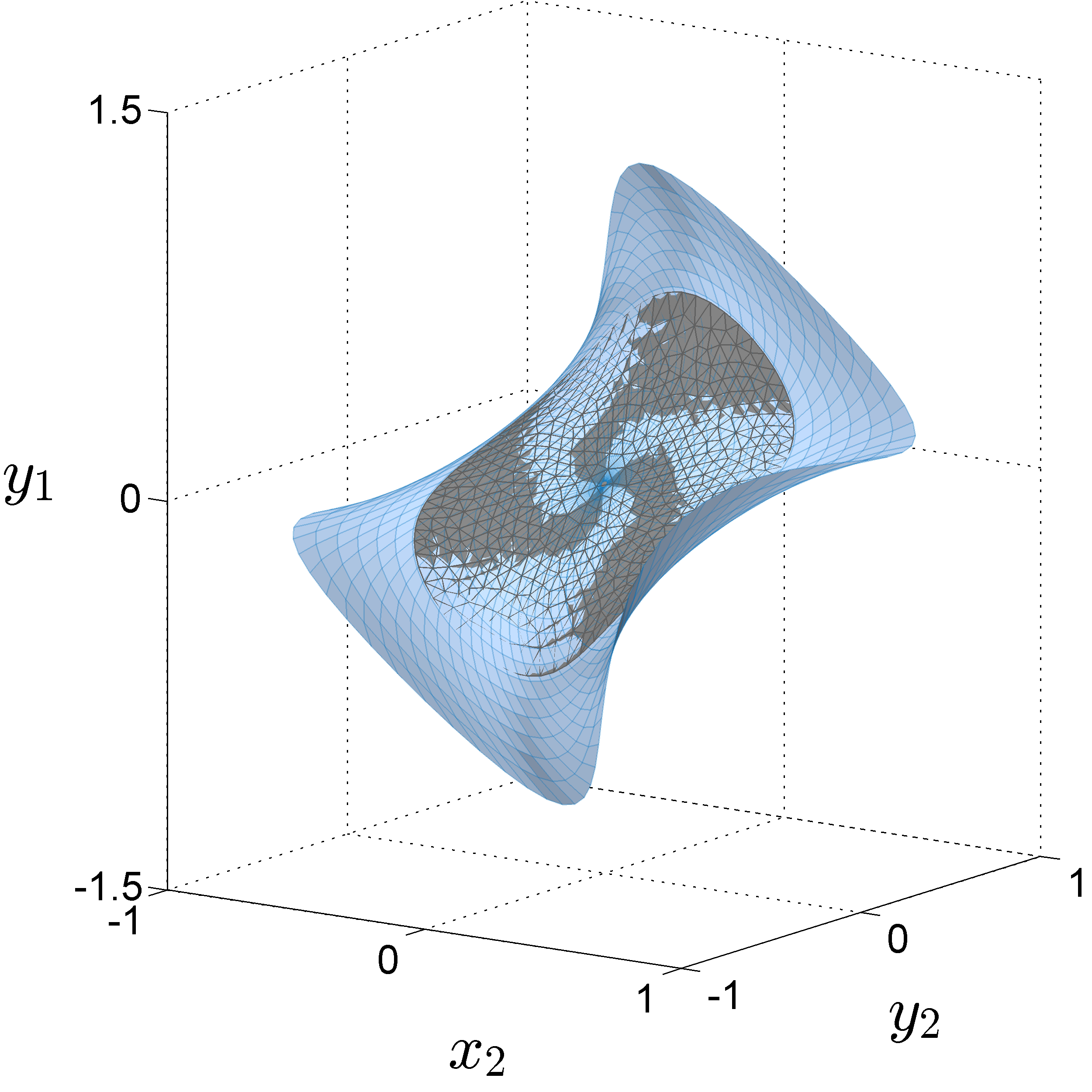}
}

\vfill

\subfloat[\label{fig:invMan41-1}]{
\includegraphics[height=.28\textheight]{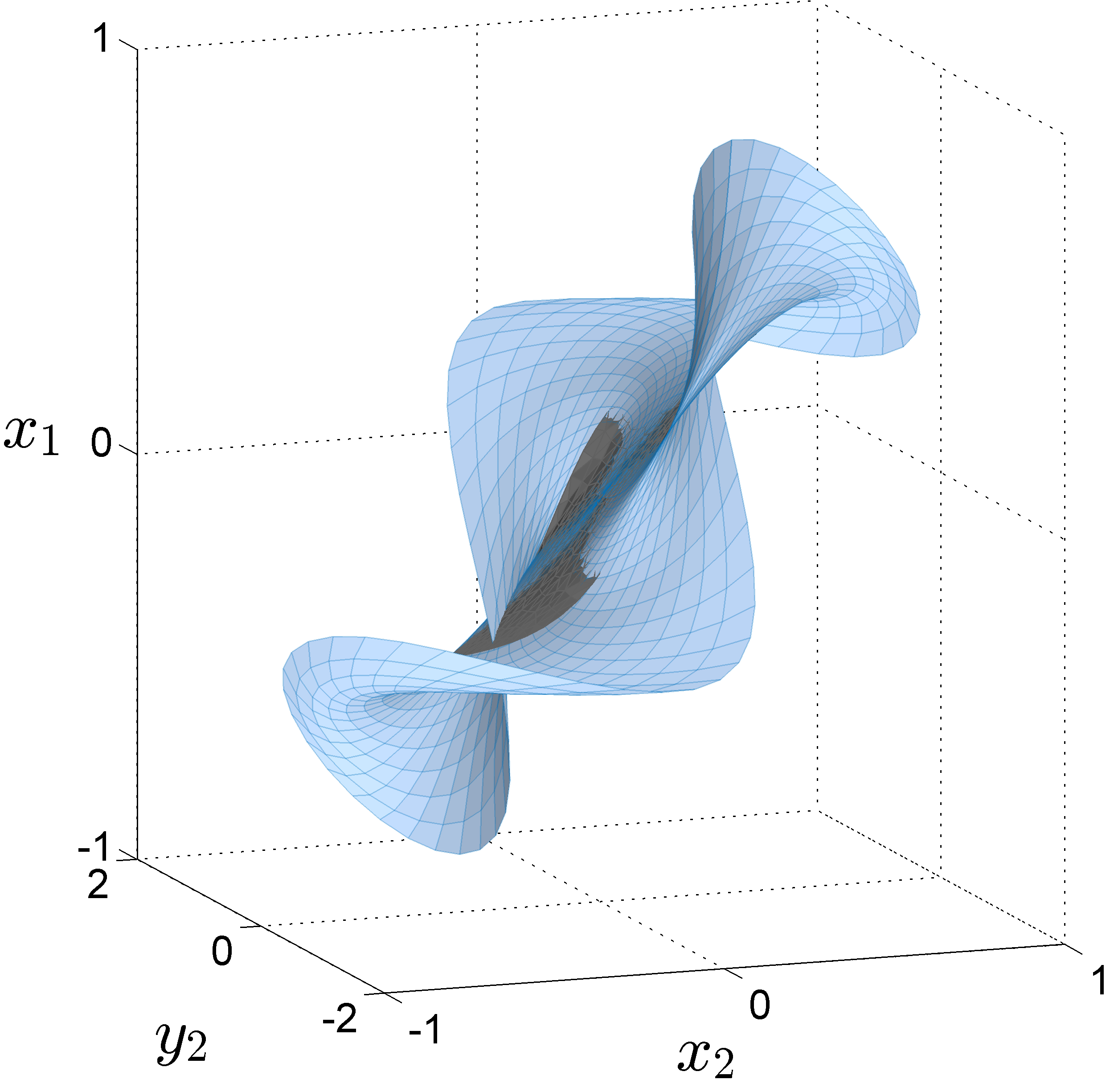}
}\hfill\subfloat[\label{fig:invMan41-2}]{
\includegraphics[height=.28\textheight]{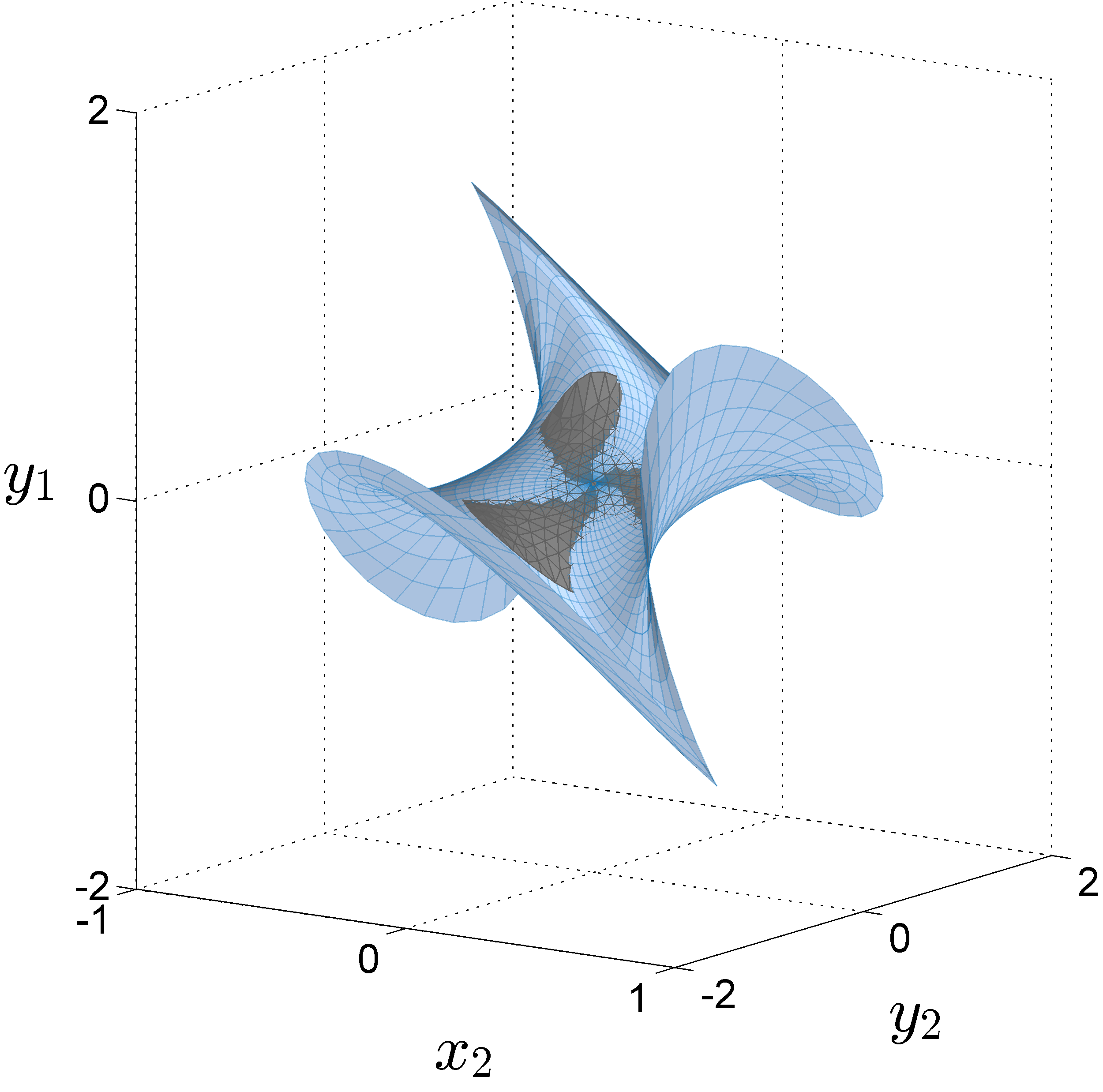}
}
\caption[In-phase NNM for different values of $k_b$]{In-phase NNM for different values of $k_b$; (a)-(b) $k_b=4.7$, (c)-(d) $k_b=4.5$, (e)-(f) $k_b=4.1$;%
\includegraphics{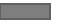}%
FE method,%
\includegraphics{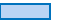}%
Taylor series}\label{fig:invManK2}
\end{figure*}

The result is compared with the finite element (FE) method developed in~\cite{renson2014effective}. This numerical method is based on an explicit parametrization of the invariant surface using a pair of master coordinates and it is known to be accurate provided that a good mesh is employed. For the out-of-phase mode, the selected master coordinates correspond to the displacement and velocity of the first mass. As can be seen in Figure~\ref{fig:firstNNM} the results from the two methods agree very well, but the FE method allows higher levels of amplitude to be reached, while the Taylor series validity is more limited.

\begin{figure}[!ht]
\centering
\subfloat[]{
\includegraphics{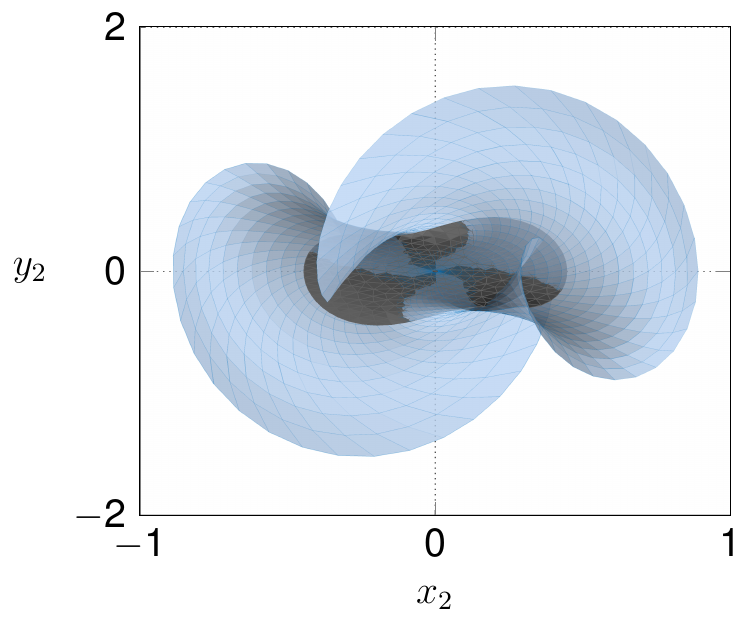}
}\subfloat[]{
\includegraphics{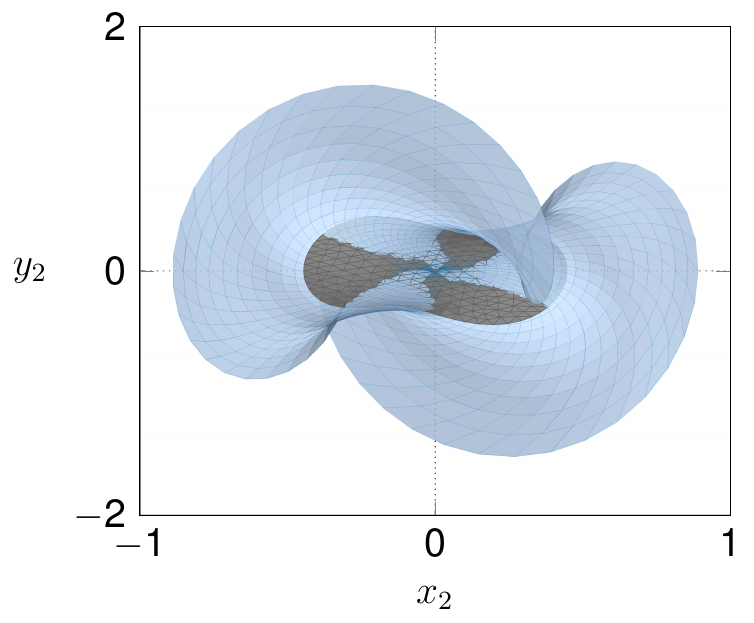}
}\caption{Different views of the invariant manifold shown in Figure~\ref{fig:invManK2} ($k_b=4.1$ ); (a) $x_1$ (Figure~\ref{fig:invMan41-1}) top view, (b) $x_2$ (Figure~\ref{fig:invMan41-2})  bottom view. Note that the manifold undergoes folding, which prevents a parametrization with state-space variables.}\label{fig:inManUV}
\end{figure}

The lower-frequency (in-phase) NNM was computed for the different values of $k_b$ reported in Table~\ref{tab:eig}. The results are shown in Figure~\ref{fig:invManK2} and compared with the NNMs obtained using the FE algorithm. As $k_b$ is decreased from $4.7$ (Figures~\ref{fig:invMan47-1} and~\ref{fig:invMan47-2}) to $4.3$ (Figures~\ref{fig:invMan43-1} and~\ref{fig:invMan43-2}) and $4.1$ (Figures~\ref{fig:invMan41-1} and~\ref{fig:invMan41-2}), the manifold geometry changes and starts to fold, which prevents the use of master coordinates for the parametrization. This is also confirmed by Figure~\ref{fig:inManUV}, where the manifold obtained for $k_b=4.1$ is shown with a viewpoint orthogonal to the $x_2y_2$ plane. The FE algorithm reaches the folding but cannot go beyond, since it cannot converge in region where the master coordinates parametrization is not valid. In contrast the folding does not pose any problem with the parametrization obtained through the Koopman eigenfunctions.

\begin{figure}[!hb]
\centering%
\hspace{.7em}\subfloat[\label{fig:x1ts}]{%
\includegraphics{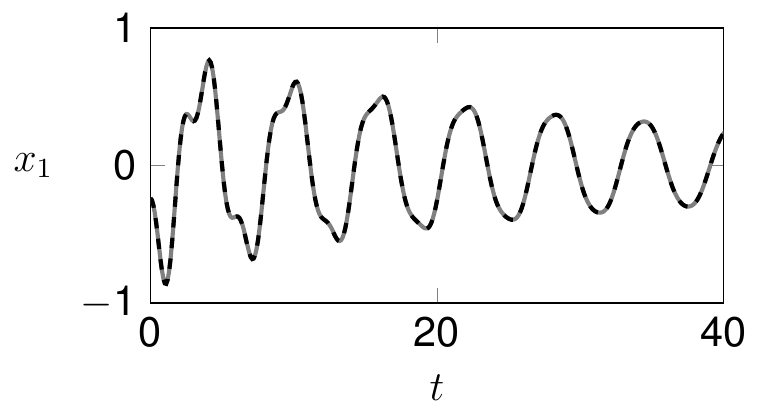}%
}\hfill\subfloat[\label{fig:x2ts}]{%
\includegraphics{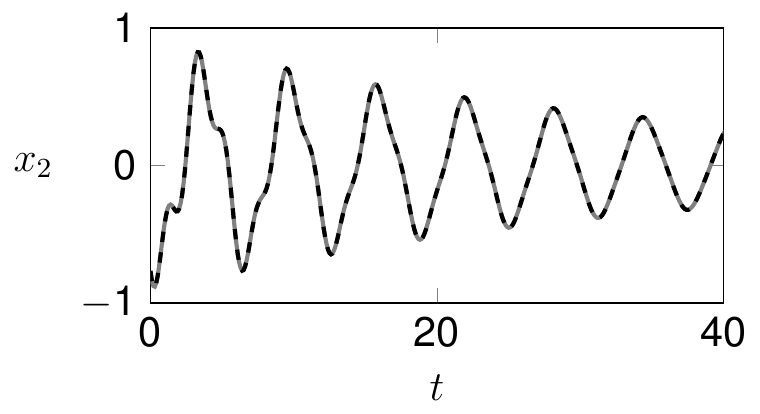}%
}

\subfloat[\label{fig:y1ts}]{%
\includegraphics{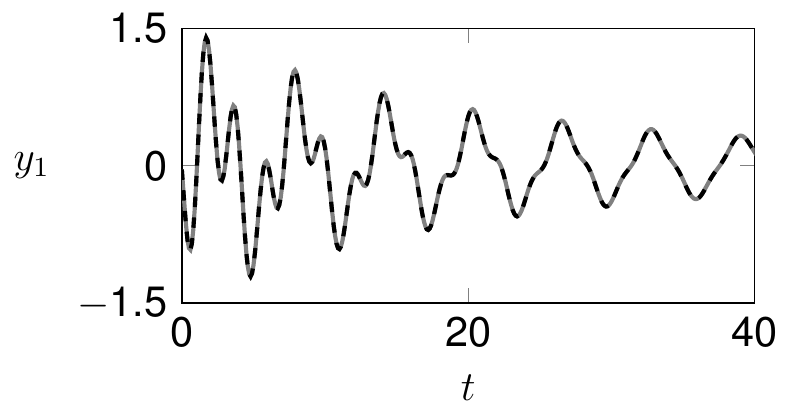}%
}\hfill\subfloat[\label{fig:y2ts}]{%
\includegraphics{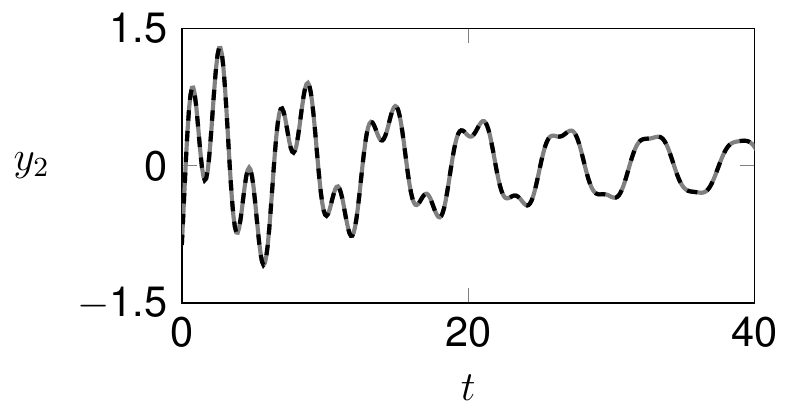}%
}%
\caption[Comparison of time series]{Comparison of time series;%
\includegraphics{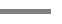}%
Numerical integration,%
\includegraphics{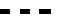}%
Taylor approximation; initial condition: $x_1(0)=-0.23$, $x_2(0)=-0.76$, $y_1(0)=-0.4$, $y_2(0)=-0.86$; (a) $x_1$, (b) $x_2$, (c) $y_1$, (d) $y_2$}\label{fig:timeSeries}
\end{figure}

In Figure~\ref{fig:timeSeries} one of the trajectories considered for the validation is compared with the corresponding numerical integration, in the case $k_b=4.1$. As expected according to the values of the NMSE, the two trajectories cannot be distinguished by visual inspection.

\begin{figure}[!h]
\centering
\subfloat[]{
\includegraphics{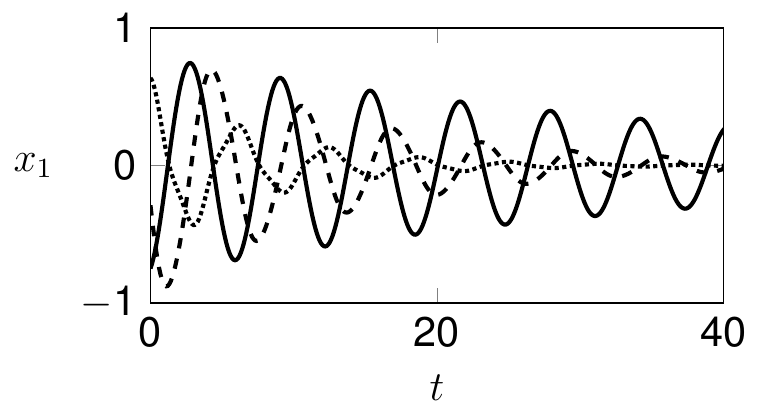}%
}\hfill\subfloat[]{
\includegraphics{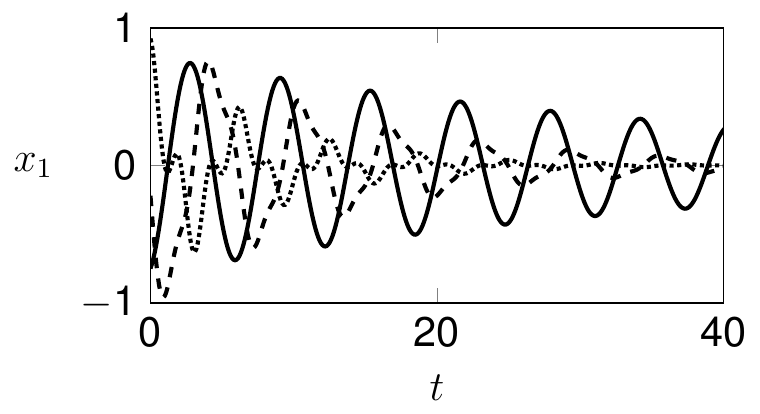}%
}
\caption[Decomposition of trajectories in different terms]{Decomposition of trajectories in different terms;%
\includegraphics{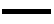}%
$1^{\text{st}}$ order term,%
\includegraphics{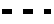}%
$3^{rd}$ order term,%
\includegraphics{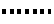}%
$5^{th}$ order term; %
(a)  $k_b=4.7$ ($x_1(0)=-0.57$, $x_2(0)=-0.55$, $y_1(0)=-0.44$, $y_2(0)=-0.66$); %
(b) $k_b=4.1$ ($x_1(0)=-0.23$, $x_2(0)=-0.76$, $y_1(0)=-0.4$, $y_2(0)=-0.86$)}\label{fig:modDec}
\end{figure}

\begin{figure}[!h]
\centering
\subfloat[]{
\includegraphics{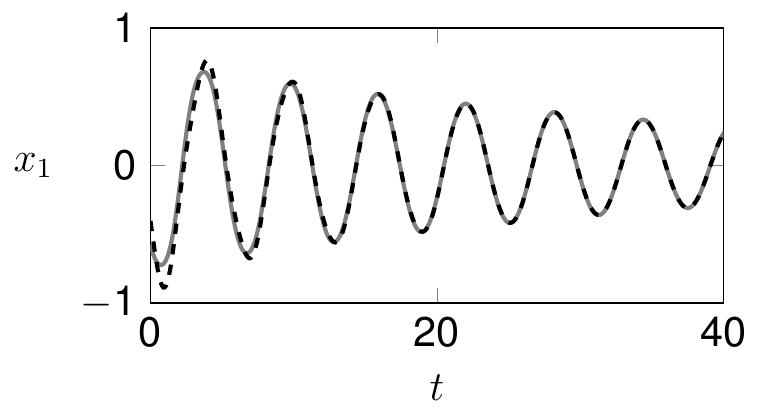}%
}\hfill\subfloat[]{
\includegraphics{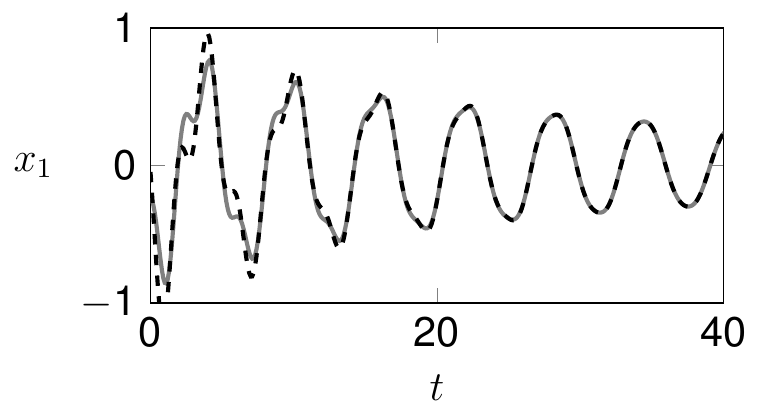}%
}
\caption[Comparison of $5^{\text{th}}$ order approximation and numerical integration]{Comparison of $5^{\text{th}}$ order approximation and numerical integration;%
\includegraphics{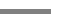}%
Numerical integration,%
\includegraphics{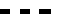}%
$5^{\text{th}}$ order approximation; %
(a) $k_b=4.7$ ($x_1(0)=-0.57$, $x_2(0)=-0.55$, $y_1(0)=-0.44$, $y_2(0)=-0.66$), %
(b) $k_b=4.1$ ($x_1(0)=-0.23$, $x_2(0)=-0.76$, $y_1(0)=-0.4$, $y_2(0)=-0.86$)}\label{fig:lowOrdApp}
\end{figure}

Equation~(\ref{eq:trajTaylor}) also highlights a structure of the evolution that is reminiscent of that obtained using Fourier series with periodic trajectories. The main difference is that the number of terms in~\eqref{eq:trajTaylor} is not the same at all orders, but increases with the order. Figure~\ref{fig:modDec} shows the terms corresponding to the first three non-zero orders, for $k_b=4.7$ and $k_b=4.1$, and in Figure~\ref{fig:lowOrdApp} their sum is compared with the trajectories obtained using numerical integration. When only terms up to order five are used the agreement between the two trajectories is not as good as in Figure~\ref{fig:timeSeries}, nonetheless, it is remarkable that the main aspects of the dynamics are retained by this approximation, at least in the case considered herein. It is also interesting to note that in the case $k_b=4.7$ the trajectory presents a lower harmonic content, and using the same order a better approximation is obtained.

\FloatBarrier%
\section{Conclusion}\label{sec:conc}
This research investigates the intimate connection that exists between NNMs of damped systems and the Koopman operator associated with their dynamics. In particular, we have shown that there is a correspondence between the invariant manifolds defined by Shaw and Pierre and the zero level sets of particular eigenfunctions of the Koopman operator, and that this correspondence is a natural generalization of the linear case. 

Through the Kooman operator eigenfunctions a new characterization of NNMs is possible, and a novel approach to their parametrization is proposed. The main advantage of this new approach is its global validity, obtained as a consequence of the linearizing properties provided by the eigenfunctions. This parametrization overcomes a known limitation in the use of master coordinates.

However, the Taylor series expansion used herein is local in nature and restricted to low-dimensional systems. In order to target arbitrarily high levels of amplitude, as well as high-dimensional systems, more efficient methods should be developed for solving the system of PDEs~\eqref{eq:NNMpdeI} or~\eqref{eq:NNMpdeR}. As mentioned in Section~\ref{sec:PDEsys} there are differences between these systems and the one obtained in the work of Shaw and Pierre. Nonetheless, all these systems share a common structure, which might allow the extension of numerical methods developed using master coordinates to this new framework. In addition, the Koopman operator approach is associated with recent data analysis techniques that are well-suited to high-dimensional systems~\cite{williams2015kernel}. This could potentially open the door to develop data-driven methods for computing NNMs of large systems.

The preliminary results presented here motivate a deeper investigation of NNMs for damped systems. For example it is well known that in conservative systems NNMs can bifurcate, and exceed the number of degrees of freedom of the system, but how these phenomena survive the presence of damping is still an open question.

\section*{Acknowledgment}
This paper presents research results of the Belgian Network DYSCO (Dynamical Systems, Control, and Optimization), funded by the Interuniversity Attraction Poles Programme initiated by the Belgian Science Policy Office.

This work was performed while A. Mauroy held a return grant from the Belgian Science Policy (BELSPO). The author L. Renson is a Marie-Curie COFUND Postdoctoral Fellow of the University of Li\`{e}ge, co-funded by the European Union.

\renewcommand{\thesection}{\Alph{section}}
\setcounter{section}{0}
\renewcommand*{\theHsection}{app.\the\value{section}}
\section{Computation of the Koopman Modes of the Identity}\label{ap:KMcomp}
In this appendix we derive a method for computing the Koopman modes $v_{k_1,\ldots,k_n}$ of the identity (see~\eqref{eq:id_dec}). In particular, we show that these modes  are the solutions of linear systems.

We fix some notation. For $\alpha=(\alpha_1,\ldots,\alpha_n)\in\mathbb{N}^n$ let 
\begin{align}
\alpha!&=\prod_{j=1}^n\alpha_j!, \\
|\alpha|&=\sum_{j=1}^n\alpha_j, \\
\partial^\alpha \Vf_l&=\frac{\partial^{|\alpha|} \Vf_l}{\partial
x_1^{\alpha_1}\ldots\partial x_n^{\alpha_n}}, \\
x^\alpha&=\prod_{j=1}^nx_j^{\alpha_j}.
\end{align}
Using this notation, we can write 
\begin{equation}
s_{k_1,\ldots,k_n}(\Vx)=s^{k}(\Vx), \,\,\,\,\,\, k=(k_1,\ldots,k_n),
\end{equation}
with
\begin{equation}
\Vs(\Vx)=\left(%
\begin{array}{c}
s_1(\Vx) \\ 
\vdots \\ 
s_n(\Vx)%
\end{array}%
\right).
\end{equation}

In addition,~\eqref{eq:id_dec} can be rewritten as 
\begin{equation}  \label{eq:id_dec_2}
\Vx=\sum_{h\geq 1}\sum_{\substack{ k\in\mathbb{N}^n \\ |k|=h}}\Vv_{k}s^k(\Vx),
\end{equation}
and 
\begin{equation}  \label{eq:F_1}
\left.\dot{\Vphi}(t,\Vx)\right|_{t=0}=\Vf(\Vx)=\sum_{h\geq 1}\sum_{\substack{ k\in%
\mathbb{N}^n \\ |k|=h}}\Vv_{k}s^k(\Vx)k\lambda,
\end{equation}
where $k\lambda=\sum_{j=1}^n k_j\lambda_j$ and $\lambda_j$ are the eigenvalues of the Jacobian matrix
\begin{equation}
\VA=\left.\frac{\partial \Vf}{\partial \Vx}\right|_{\Vx=\Vzero}.
\end{equation}
Similarly, the Taylor series of $\Vf(\Vx)$ can be written as 
\begin{equation}\label{eq:Ftaylor}
\Vf_l(\Vx)=\sum_{h\geq 1}\sum_{\substack{ k\in\mathbb{N}^n \\ |k|=h}}\frac{1}{k!}%
\partial^k \Vf_l(\Vzero)\Vx^k,
\end{equation}
and injecting~\eqref{eq:id_dec_2} in~\eqref{eq:Ftaylor} leads to 
\begin{equation}  \label{eq:F_l2}
\Vf_l(\Vx)=\sum_{h_1\geq 1}\sum_{\substack{ k_1\in\mathbb{N}^n \\ |k_1|=h_1}}%
\frac{1}{k_1!}\partial^{k_1} \Vf_l(\Vzero)\prod_{j=1}^n\left(\sum_{h_2>1}\sum_{\substack{
k_2\in\mathbb{N}^n \\ |k_2|=h_2}}(\Vv_{k_2})_j s^{k_2}(\Vx)\right)^{k_{1_j}}.
\end{equation}

By comparison of~\eqref{eq:F_1} and~\eqref{eq:F_l2}, for any $k\in \mathbb{N}^{n}$ we obtain 
\begin{equation}\label{eq:idModeLinSyst}
(\Vv_{k})_{l}k\lambda =\sum_{h=1}^{|k|}\sum_{|j|=h}\frac{1}{j!}\partial
^{j}\Vf_{l}(\Vzero)\sum_{\substack{ (k_{1}\ldots k_{|j|})| \\ \sum_{i}k_{i}=k}}%
\prod_{i=1}^{n}\prod_{m=\sum\limits_{s}^{i-1}j_{s}+1}^{\sum%
\limits_{s}^{i}j_{s}}(\Vv_{k_{m}})_{i}.
\end{equation}%
These equations are linear in the terms with maximum value of $|k|$, and can
be solved if all the $\Vv_{k^{\prime }}$ with $|k^{\prime }|<|k|$ are known. 

The matrix defining the linear system related to $\Vv_k$ has the form $k\lambda \VI-\VA$. The first-order modes are the eigenvectors of $\VA$ and it is necessary to choose a scaling factor, which correspondingly fixes the value of the eigenfunctions (defined up to a multiplicative constant). If there are no resonances between the eigenvalues, the  linear systems related to modes of higher order are uniquely determined, since they are defined by a nonsingular matrix. We can also use the known form of the matrices defining the linear systems to simplify the computation. In the case of a diagonalizable Jacobian matrix, the systems are diagonal when expressed in coordinates defined by a Jordan basis.

To obtain a manifold using~\eqref{eq:NNM_taylor} not all the modes are necessary, but just a subset of them: the $\Vv_k$, $k=(k_1,\ldots,k_n)$, for which only two indices $k_{j_1}$ and $k_{j_2}$ are nonzero. Since to obtain the linear system corresponding to a given mode $\Vv_k$ only the modes $\Vv_{k'}$, $k'=(k'_1,\ldots,k'_n)$, with $k_l\ge k'_l$ are needed, the computation can be restricted to the modes that appear in~\eqref{eq:NNM_taylor}.

As concluding remarks we note that approximating a solution of~\eqref{eq:NNMpdeI} with a Taylor series yields the same linear systems. Also, the procedure presented here corresponds in essence to a normal form approximation. In that context, starting from a dynamical system with diagonal linear part, a transformation can be obtained so that the resulting equations of motion are the simplest possible. The only algebraic condition needed to obtain a transformation that linearizes the dynamic is the nonresonant eigenvalues. The corresponding linearizing coordinates are the eigenfunctions of the Koopman operator.

\section{A family of invariant manifolds for a linear system}\label{ap:mm}
Consider a four-dimensional linear system 
\begin{equation}\label{eq:vecField}
\dot{\Vx}=\VA\Vx,
\end{equation}
with
\begin{equation}
\VA=\left(
\begin{array}{cccc}
-\sigma_1 &  \omega_1 & 0 & 0 \\
-\omega_1 & -\sigma_1 & 0 & 0 \\
0 & 0 & -\sigma_2 &  \omega_2 \\
0 & 0 & -\omega_2 & -\sigma_2
\end{array}\right).
\end{equation}
Note that any four-dimensional linear system with two couples of complex conjugate eigenvalues can be transformed into this form. It will be assumed that $\sigma_1>0$ and $\sigma_2>0$ (stability).

For any value of $c_1$ and $c_2$ consider the manifold $\Vx(u,v)=\Vf(u,v;c_1,c_2)$,
\begin{align}\label{eq:invLinMan}
\begin{split}
f_1(u,v)&= u \\
f_2(u,v)&= v \\
f_3(u,v)&= \left(u^2+v^2\right)^{\frac{\sigma _2}{2 \sigma _1}} \left(\cos
   \left(\frac{\log \left(u^2+v^2\right) \omega _2}{2 \sigma
   _1}\right) c_1-\sin \left(\frac{\log \left(u^2+v^2\right)
   \omega _2}{2 \sigma _1}\right) c_2\right) \\
f_4(u,v)&= \left(u^2+v^2\right)^{\frac{\sigma _2}{2 \sigma _1}} \left(\sin
   \left(\frac{\log \left(u^2+v^2\right) \omega _2}{2 \sigma
   _1}\right) c_1+\cos \left(\frac{\log \left(u^2+v^2\right)
   \omega _2}{2 \sigma _1}\right) c_2\right).
\end{split}
\end{align}

The derivatives of $\Vf$ are
\begin{equation}
\begin{split}
f_{u,1}&=  1 \\
f_{u,2}&= 0 \\
f_{u,3}&= \frac{u \left(u^2+v^2\right)^{\frac{\sigma _2}{2 \sigma _1}-1}
   \left(\cos \left(\frac{\log \left(u^2+v^2\right) \omega _2}{2
   \sigma _1}\right) \left(c_1 \sigma _2-c_2 \omega _2\right)-\sin
   \left(\frac{\log \left(u^2+v^2\right) \omega _2}{2 \sigma
   _1}\right) \left(c_2 \sigma _2+c_1 \omega
   _2\right)\right)}{\sigma _1} \\
f_{u,4}&= \frac{u \left(u^2+v^2\right)^{\frac{\sigma _2}{2 \sigma _1}-1}
   \left(\cos \left(\frac{\log \left(u^2+v^2\right) \omega _2}{2
   \sigma _1}\right) \left(c_2 \sigma _2+c_1 \omega _2\right)+\sin
   \left(\frac{\log \left(u^2+v^2\right) \omega _2}{2 \sigma
   _1}\right) \left(c_1 \sigma _2-c_2 \omega
   _2\right)\right)}{\sigma _1} \\
f_{v,1}&= 0 \\
f_{v,2}&= 1 \\
f_{v,3}&= \frac{v \left(u^2+v^2\right)^{\frac{\sigma _2}{2 \sigma _1}-1}
   \left(\cos \left(\frac{\log \left(u^2+v^2\right) \omega _2}{2
   \sigma _1}\right) \left(c_1 \sigma _2-c_2 \omega _2\right)-\sin
   \left(\frac{\log \left(u^2+v^2\right) \omega _2}{2 \sigma
   _1}\right) \left(c_2 \sigma _2+c_1 \omega
   _2\right)\right)}{\sigma _1} \\
f_{v,4}&= \frac{v \left(u^2+v^2\right)^{\frac{\sigma _2}{2 \sigma _1}-1}
   \left(\cos \left(\frac{\log \left(u^2+v^2\right) \omega _2}{2
   \sigma _1}\right) \left(c_2 \sigma _2+c_1 \omega _2\right)+\sin
   \left(\frac{\log \left(u^2+v^2\right) \omega _2}{2 \sigma
   _1}\right) \left(c_1 \sigma _2-c_2 \omega
   _2\right)\right)}{\sigma _1}.
\end{split}
\end{equation}
These are continuous for $\sigma_2>\sigma_1$, and more regular if the distance between the two real parts of the eigenvalues is increased further.

Assuming sufficient regularity, $\Vf$ defines a manifold, with tangent plane at a point $(u,v)$ spanned by the vectors
\begin{equation}
\Vf_u=\frac{\partial \Vf}{\partial u},\,\,\,\,\Vf_v=\frac{\partial \Vf}{\partial v}.
\end{equation}
Since it holds
\begin{equation}\label{eq:invEq}
\VA\Vf(u,v)=(-\sigma_1 u +\omega_1 v)\Vf_u+(-\omega_1 u -\sigma_1 v)\Vf_v,
\end{equation}
all the manifolds are invariant (\eqref{eq:invEq} is the analogous of the equation derived by Shaw and Pierre equation for the system~\eqref{eq:vecField}). Moreover, if $\sigma_2>\sigma_1$
\begin{equation}
\Vf_u(0,0)=\left(\begin{array}{c}
1\\0\\0\\0
\end{array}\right),\,\,\,\, \Vf_v(0,0)=\left(\begin{array}{c}
0\\1\\0\\0
\end{array}\right),
\end{equation}
so that at the origin the tangent plane is $x_3=0,x_4=0$, which corresponds to the linear normal mode, but only for $c_1=c_2=0$ the whole manifold is the linear normal mode.

Note that apart from the points verifying $x_{1,0}=x_{2,0}=0$ (i.e. points on the fast linear mode), for any other point $(x_{1,0},x_{2,0},x_{3,0},x_{4,0})$, there exists a manifold which passes through it, this is obtained using the parameters
\begin{equation}
\begin{split}
c_1&=x_{3,0} \left(\frac{1}{x_{1,0}^2+x_{2,0}^2}\right){}^{\frac{\sigma
   _2}{2 \sigma _1}} \cos \left(\frac{\omega _2 \log
   \left(\frac{1}{x_{1,0}^2+x_{2,0}^2}\right)}{2 \sigma
   _1}\right)-x_{4,0}
   \left(\frac{1}{x_{1,0}^2+x_{2,0}^2}\right){}^{\frac{\sigma
   _2}{2 \sigma _1}} \sin \left(\frac{\omega _2 \log
   \left(\frac{1}{x_{1,0}^2+x_{2,0}^2}\right)}{2 \sigma _1}\right),\\
c_2&=x_{3,0} \left(\frac{1}{x_{1,0}^2+x_{2,0}^2}\right){}^{\frac{\sigma
   _2}{2 \sigma _1}} \sin \left(\frac{\omega _2 \log
   \left(\frac{1}{x_{1,0}^2+x_{2,0}^2}\right)}{2 \sigma
   _1}\right)+x_{4,0}
   \left(\frac{1}{x_{1,0}^2+x_{2,0}^2}\right){}^{\frac{\sigma
   _2}{2 \sigma _1}} \cos \left(\frac{\omega _2 \log
   \left(\frac{1}{x_{1,0}^2+x_{2,0}^2}\right)}{2 \sigma _1}\right).
\end{split}
\end{equation}

As a final remark we note that even if each couple of coefficients $c_1$, $c_2$ defines a different manifold, the derivatives of all the functions defined by~\eqref{eq:invLinMan} have the same values at the origin (up to the order at which they are defined): they all share the same polynomial approximation. Interestingly the same situation is encountered in the context of center manifolds.


\bibliographystyle{model1-num-names}
\bibliography{bibKONNM}

\end{document}